\documentclass[11pt]{article}
\usepackage{graphicx}
\usepackage{color}
\usepackage{amsmath}
\usepackage{amssymb}
\usepackage{amscd}
\usepackage{bbm}
\usepackage{framed}

\usepackage{fancyhdr,a4wide}
\usepackage{amsthm}

\newcommand{\R}{\mathbb{R}}
\newcommand{\inr}[1]{\left\langle #1 \right\rangle}

\newcommand{\E}{\mathbb{E}}
\newcommand{\M}{\mathbb{M}}
\newcommand{\Q}{\mathbb{Q}}

\newcommand{\eps}{\varepsilon}

\newtheorem{Theorem}{Theorem}[section]
\newtheorem{Lemma}[Theorem]{Lemma}

\newtheorem{Definition}[Theorem]{Definition}

\newtheorem{Proposition}[Theorem]{Proposition}
\newtheorem{Corollary}[Theorem]{Corollary}
\newtheorem{Remark}[Theorem]{Remark}

\newtheorem{Assumption}{Assumption}[section]

\numberwithin{equation}{section}

\def \endproof
{{\mbox{}\nolinebreak\hfill\rule{2mm}{2mm}\par\medbreak}}

\newcommand{\ol}{\overline}
\newcommand{\wt}{\widetilde}
\newcommand{\wh}{\widehat}
\newcommand{\argmin}{\mathop{\mathrm{argmin}}}
\newcommand{\X}{\mathcal{X}}
\newcommand{\C}{\mathcal{C}}
\newcommand{\D}{\mathcal{D}}

\newcommand{\EXP}{\mathbb{E}}
\newcommand{\PROB}{\mathbb{P}}

\newcommand{\F}{{\mathcal F}}
\newcommand{\cH}{{\mathcal H}}

\begin{document}
\title{Regularization, sparse recovery, and median-of-means tournaments
\thanks{
G\'abor Lugosi was supported by
the Spanish Ministry of Economy and Competitiveness,
Grant MTM2015-67304-P and FEDER, EU. Shahar Mendelson was supported in part by the Israel Science Foundation.
}
}
\author{
G\'abor Lugosi\thanks{Department of Economics and Business, Pompeu
  Fabra University, Barcelona, Spain, gabor.lugosi@upf.edu}
\thanks{ICREA, Pg. Lluís Companys 23, 08010 Barcelona, Spain}
\thanks{Barcelona Graduate School of Economics}
\and
Shahar Mendelson \thanks{Department of Mathematics, Technion, I.I.T, and Mathematical Sciences Institute, The Australian National University, shahar@tx.technion.ac.il}}

\maketitle

\begin{abstract}
We introduce a regularized risk minimization procedure for regression function estimation. The procedure is based on median-of-means tournaments, introduced by the authors in
\cite{LuMe16} and achieves near optimal accuracy and confidence under general conditions,
including heavy-tailed predictor and response variables. It outperforms standard regularized
empirical risk minimization procedures such as {\sc lasso} or {\sc slope} in heavy-tailed problems.
\end{abstract}


\section{Introduction}

\subsection{Empirical risk minimization, regularization}

Regression function estimation is a fundamental problem in statistics and machine learning.
In the most standard formulation of the problem, $(X,Y)$ is a pair of random variables in which
$X$, taking values in some general measurable space $\X$, represents the observation (or feature vector)
and one would like to approximate the unknown real value $Y$ by a function of $X$.
In other words, one is interested in finding a function $f:\X\to \R$ such that $f(X)$ is ``close'' to $Y$.
As the vast majority of the literature, we measure the quality of $f$ by the \emph{risk}
\[
    R(f) = \EXP (f(X)-Y)^2~,
\]
which is well defined whenever $f(X)$ and $Y$ are square integrable, assumed throughout the paper.
Clearly, the best possible function is the \emph{regression function} $m(X)= \EXP(Y|X)$.

However, in statistical problems, the joint distribution of $(X,Y)$ is unknown and the regression function
is impossible to compute. Instead, a sample $\D_N=((X_1,Y_1),\ldots,(X_N,Y_N))$ of independent copies
of the pair $(X,Y)$ is available (such that $\D_N$ and the pair
$(X,Y)$ are independent).

A popular and thoroughly studied approach is to select a function $\wh{f}_N$ from a fixed class $\F$ of functions.
Formally, a \emph{learning procedure} is a map
$\Phi:(\X \times \R)^N \to \F$ that assigns to each sample
$\D_N=(X_i,Y_i)_{i=1}^N$ a (random) function $\Phi(\D_N)=\wh{f}_N$. If the class $\F$ is sufficiently ``large,'' then it is reasonable to expect that the best function in the class
\[
      f^* = \argmin_{f\in\F} \EXP (f(X)-Y)^2
\]
has an acceptable performance, and the standard assumption is that the minimum is attained and $f^*\in \F$ is unique. We assume that $\F$ is a closed and convex subset of $L_2(\mu)$---where $\mu$ denotes
the distribution of $X$---, guaranteeing the existence and uniqueness of $f^*$.

The quality of a learning procedure is typically measured by the \emph{mean squared error}, which is the conditional expectation
$$
\|\wh{f}_N-f^*\|_{L_2}^2= \E \bigl((\wh{f}_N(X)-f^*(X))^2 | \D_N \bigr) ~,
$$
where, for $q\ge 1$, we use the notation
\[
 \|f-g\|_{L_q}=    \left(\EXP \left|f(X)-g(X)\right|^q \right)^{1/q}
\quad \text{and also} \quad
 \|f-Y\|_{L_q}=    \left(\EXP \left|f(X)-Y\right|^q \right)^{1/q}~.
\]
A closely related, though not equivalent, measure of performance is the \emph{excess risk}, defined by the conditional expectation
$$
R(\wh{f}_N) - R(f^*) =\E \bigl((\wh{f}_N(X)-Y)^2|\D_N\bigr) - \E (f^*(X)-Y)^2~.
$$
The goal of a statistical learning problem is to find a learning procedure that achieves a good \emph{accuracy}
with a high \emph{confidence}. In particular, for $r>0$ and $\delta \in (0,1)$,
we say that a procedure performs with accuracy parameter $r$ with confidence $1-\delta$ in the class $\F$ (e.g., for the mean squared error)
if
\[
      \PROB \left( R(\wh{f}_N)-R(f^*) \le r^2 \right) \ge 1- \delta~
\]
(sometimes one only considers $\PROB \left( \|\wh{f}_N-f^*\|_{L_2} \le r \right) \ge 1- \delta$).

High accuracy and high confidence (i.e., small $r$ and small $\delta$) in the given class are obviously conflicting requirements.
The achievable tradeoff has been thoroughly studied and it is fairly well understood. We refer the reader to
Lecu\'e and Mendelson \cite{LeMe16}, Lugosi and Mendelson \cite{LuMe16} for recent accounts.

The most standard approach for a learning procedure is \emph{empirical risk minimization} ({\sc erm}), also
known as \emph{least squares regression} in which
\[
  \wh{f}_N \in \argmin_{f\in \F} \sum_{i=1}^N (f(X_i)-Y_i)^2
\]
 (where we assume that the minimum is achieved). One may show (see, e.g., Lecu\'e and Mendelson \cite{LeMe16}) that unless the function class and target are sub-Gaussian\footnote{Here, Sub-Gaussian means that the $\psi_2$ and the $L_2$ norms are equivalent in $\F \cup \{0\}$; that is, there is a constant $L$ such that for any $f,h \in \F \cup \{0\}$ and any $p \geq 2$, $\|f-h\|_{L_p} \leq L\sqrt{p}\|f-h\|_{L_2}$, and that the same holds for any $Y-f(X)$ for any $f \in \F$.}, then empirical risk minimization is far from achieving the optimal accuracy/confidence tradeoff. The reason for the suboptimal behaviour of empirical risk minimization is that outliers distort the empirical means unless the problem is very close to being Gaussian. Thankfully, learning procedures that can tackle heavy-tailed problems exist, as it was recently pointed out by Lugosi and Mendelson \cite{LuMe16} with the introduction of the \emph{median-of-means tournament}.

A common problem that all learning procedures encounter is that of overfitting, which occurs when the underlying class is too big relative to the (random) information at the learner's disposal. A standard way of dealing with learning problems involving classes that are too large is giving priority to functions in the class according to some prior belief of ``simplicity''. For example, in regularized risk minimization, one selects a norm $\Psi$ defined on a vector space $E$ containing $\F$. A small value of $\Psi(f)$ is
interpreted as simplicity and simple functions are given priority by way of adding a penalty term to the empirical risk
that is proportional to $\Psi(f)$. In particular, for some \emph{regularization parameter} $\lambda >0$,
a regularized risk minimizer selects
\[
  \wh{f}_N \in \argmin_{f\in \F} \left( \frac{1}{N} \sum_{i=1}^N (f(X_i)-Y_i)^2 + \lambda \Psi(f) \right)~,
\]
and the term $\Psi(f)$ is sometimes called the \emph{penalty}.

Just as the tournament procedure from \cite{LuMe16} outperforms
empirical risk minimization (in fact, the tournament procedure attains
the optimal tradeoff between accuracy and confidence under minimal
assumptions), the regularized tournament which we present here,
outperforms regularized risk minimization. Since regularized
procedures require the minimization of a functional that has the
empirical mean as a component, they suffer from the same disadvantages
as empirical risk minimization. Therefore, the accuracy/confidence
tradeoff exhibited by regularized risk minimization is suboptimal once
one leaves the sub-Gaussian realm, and deteriorates further if the
problem is more heavy-tailed. In contrast, we show that the
regularized tournament attains the optimal accuracy/confidence
tradeoff under rather minimal conditions and, in particular, in
heavy-tailed problems.

The paper is organized as follows.
In Section \ref{sec:procedure} we introduce a new regularized ``tournament''
procedure in a quite general framework and illustrate how it works
on an important specific case, the tournament {\sc lasso}
(see Section \ref{sec-lasso-intro}).
In Section \ref{sec:main} the main general performance bound is presented
for the regularized tournament procedure under certain specific choice
of the parameters of the procedure.
The proof of the main result is detailed in Section \ref{sec:proof}.
Finally, in Section \ref{sec:examples} two examples are worked out.
The first is a ``tournament'' version of {\sc lasso}
(introduced in Section \ref{sec-lasso-intro}) and the second is
the tournament {\sc slope}, a generalized version of tournament {\sc lasso}.

\section{The procedure}
\label{sec:procedure}

Let us now describe the regularized tournament procedure. Recall that the learner is given a closed and convex class of functions $\F \subset L_2(\mu)$ and a regularization function $\Psi$ which is assumed to be a norm on ${\rm span}(\F)$.

An underlying assumption is that $\F$ can be naturally decomposed to a hierarchy of subclasses, that is, there is a finite decreasing collection of subsets of $\F$, denoted by $(\F_\ell)_{\ell=1}^K$, such that
\[
\F = \F_1 \supset\cdots \supset \F_K~.
\]
The idea is that the sets $\F_\ell$ capture some notion of
`complexity': the larger $\ell$ is, the simpler the functions in
$\F_\ell$ are.

An important example of a hierarchy, discussed below in detail,
is based on the notion
of sparsity in $\R^d$ relative to a fixed orthonormal basis
$(e_i)_{i=1}^d$: the vector $t \in \R^d$ is $s$-sparse if its
representation in the basis $(e_i)_{i=1}^d$, $\sum_{i=1}^d t_ie_i$,
has at most $s$ nonzero coefficients.  For the class of linear
functionals $\F = \{\inr{t,\cdot} : t \in \R^d\}$ one may set $\F_\ell
=\{\inr{t,\cdot} : t \ {\rm is } \ d/2^{\ell-1} \ {\rm sparse} \}$.

Given a hierarchy $(\F_\ell)_{\ell=1}^K$, the first, second and third
phases of the procedure are performed on each one of the subclasses
$\F_\ell$ separately and the procedure returns subclasses $\cH_\ell
\subset \F_\ell$ consisting of functions whose statistical performance
is good enough to be considered as candidates for selection. The
fourth and final phase compares the candidates selected from each
class $\F_\ell$ and selects one of them, as we describe below.

The first three stages of the procedure use independent data.
In order to accommodate this, one needs to split the available data
into three independent parts. For simplicity of the presentation, we
assume that these parts have equal size, each containing $N$ samples.
(Thus, the total sample size is $3N$ rather than $N$ but this change
of convention only affects the constants in the bounds that we do not
make explicit in any case.)

\subsubsection*{The first phase: constructing the `referee' in $\F_\ell$}

The goal of the first phase is to get a data-dependent estimate of
$L_2(\mu)$ distances between elements in $\F_\ell$. Its output is a
(data-dependent) function ${\cal DO}_\ell : \F_\ell \times \F_\ell \to
\{0,1\}$, which takes the value $1$ when two functions in $\F_\ell$
are far-enough, in the following sense. This first step uses unlabeled
samples only.
\begin{framed}
\begin{Definition} \label{def:first-phase}
  Suppose we are given a sample $(X_i)_{i=1}^{N}$.
  Fix an appropriately chosen positive integer $n_1$.
For each $\ell=1,\ldots,K$,
\begin{description}
\item{$(1)$} fix well-chosen parameters $\rho_\ell$ and $r_{\ell,1}$;

\item{$(2)$} split $(X_i)_{i=1}^N$ to $n_1$ disjoint blocks $(I_j)$ of equal size, denoted by $m_1=N/n_1$;
\item{$(3)$} for every $f,h \in \F_\ell$ let
$$
v_j= \frac{1}{m_1} \sum_{i \in I_j} |f(X_i)-h(X_i)|~,
$$
and set $d_\ell(f,h)$ to be the median of $\{ v_j : 1 \leq j \leq n_1\}$;
\item{$(4)$} set ${\cal DO}_\ell(f,h)=1$ if either $\Psi(f-h) \geq \rho_\ell$ or if $\Psi(f-h) < \rho_\ell$ and $d_\ell(f,h) \geq r_{\ell,1}$.
\end{description}
\end{Definition}
\end{framed}

\subsubsection*{The second phase: $\ell$-elimination}
This first phase is used as input to the second phase. The latter is a
comparison of the statistical performance of any two functions $f$ and
$h$ in $\F_\ell$ using a sample $(X_i,Y_i)_{i=N+1}^{2N}$. The idea is
to define `statistical matches' between any two functions $f,h \in
\F_\ell$, with each match designed to determine which function out of
the two is more suited for our needs. Because the `matches' are based
on a random sample, a reliable comparison is impossible when the two
functions are `too close'. This is where the output of the first phase
comes into the frame. The random binary function ${\cal DO}_\ell$
tells the learner when $f$ and $h$ are too close (when ${\cal
  DO}_\ell(f,h)=0$) and in such cases, the outcome of the statistical
match between $f$ and $h$ is useless.

More accurately, the second phase is defined as follows:
\begin{framed}
\begin{Definition} \label{def:phase-2}
Suppose we are given a sample $(X_i,Y_i)_{i=N+1}^{2N}$ and, for each $\ell=1,\ldots,K$, the output ${\cal DO}_\ell$ of the first phase, Then
\begin{description}
\item{$(1)$} let $\lambda_\ell$ and $n_{\ell,2}$ be well-chosen parameters, and set $\rho_\ell$ as in the first phase;
\item{$(2)$} split $\{N+1,\ldots,2N\}$ to $n_{\ell,2}$ coordinate blocks $(I_j)$ of equal size, denoted by $m_{\ell,2}$;
\item{$(3)$} for any $f,h \in \F_\ell$ set $f \succ h$ if
$$
\frac{1}{m_{\ell,2}} \sum_{i \in I_j} \left((h(X_i)-Y_i)^2 - (f(X_i)-Y_i)^2\right) + \lambda_\ell \left(\Psi(h)-\Psi(f)\right) > 0
$$
for the majority of the blocks $I_j$;
\item{$(4)$} denote by
$$
\cH_\ell^\prime = \{ f \in \F_\ell \ : \ f \succ h \ {\rm for \ every \ } h \in \F_\ell \ {\rm such \ that } \ {\cal DO}_\ell(f,h)=1\}.
$$
\end{description}
\end{Definition}
\end{framed}

The output of the second phase is the sequence of classes
$\cH_\ell^\prime \subset \F_\ell$. Each $\cH_\ell^\prime$
consists of the set of functions
that were superior to all the other `competitors' in the
$\ell$-tournament, at least in `matches' that were considered reliable
by ${\cal DO}_\ell$. We will show that if $f^* \in \F_\ell$ then functions
in $\cH_\ell^\prime$ are `close' in $L_2(\mu)$ to $f^*$.

\subsubsection*{The third phase: $\ell$-champions league}
The third phase in the procedure is a further selection process, this time conducted among the set of `winners' of the $\ell$-elimination tournament.
\begin{framed}
  \begin{Definition}
\label{def:phase-3}
    Suppose we are given $(X_i,Y_i)_{i=2N+1}^{3N}$. For each $\ell=1,\ldots,K$,
let $\cH_\ell^\prime$ be the outcome of the second phase. Then
\begin{description}
\item{$(1)$} let $n_{\ell,2}$ be as above, and set $r_{\ell,3}$ to be a well-chosen parameter;
\item{$(2)$} split $\{2N+1,\ldots,3N\}$ to $n_{\ell,2}$ disjoint blocks $(I_j)$ of size $m_{\ell,2}=N/n_{\ell,2}$;
\item{$(3)$} for any $f,h \in \cH_\ell^\prime$ set $f \gg h$ if
$$
\frac{2}{m_{\ell,2}} \sum_{i \in I_j} (h(X_i)-f(X_i)) \cdot (f(X_i)-Y_i) \geq -r^2_{\ell,3}
$$
(observe that it is possible that $f \gg h$ and $h \gg f$ at the same time);
\item{(4)} let
$$
\cH_\ell = \{f \in \cH_\ell^\prime \ : \ f \gg h \ {\rm for \ every \ } h \in \cH_\ell^\prime \}~.
$$
\end{description}
\end{Definition}
\end{framed}
This further selection process improves the outcome of the $\ell$-elimination phase: instead of just ensuring that the functions are close to $f^*$, the $\ell$-champions league selects functions whose risk is close to the risk of $f^*$.

\subsubsection*{The fourth phase: naming a winner}
The outcome of the first three phases results in data-dependent choices $\cH_\ell \subset \F_\ell$ (and obviously at this point there is no way of knowing that the sets $\cH_\ell$ are nonempty). All the functions in each one of the $\cH_\ell$'s are in some sense the best estimators one can find within $\F_\ell$, taking into account its size and the data at the learner's disposal. Naturally, the larger $\F_\ell$ is, the larger the error one will incur by selecting an estimator in it. Therefore, in the final stage the aim is to find the smallest class in the hierarchy in which a good estimator still exists.
\begin{framed}
\begin{Definition}
Given the classes $(\cH_\ell)_{\ell=1}^K$, let $\ell_1$ be the largest integer $\ell$ such that $\bigcap_{j \leq \ell} \cH_j \not= \emptyset$. Select $\wh{f}$ to be any function in $\bigcap_{j \leq \ell_1} \cH_j$.
\end{Definition}
\end{framed}

The first three phases are the key components of the procedure. Out of the three, the first one is an adaptation of the \emph{distance oracle} used in \cite{LuMe16} and which had been introduced in \cite{Men16a}; the third component is essentially the same as the champions league stage in the tournament procedure from \cite{LuMe16}.

The truly new component is the second phase. Its analysis combines ideas from \cite{LeMe16a}  (which focused on regularized risk minimization in `sparse' problems) and from \cite{LuMe16}. As it is the main novelty in this article we present it in detail and only sketch the arguments needed in the analysis of the other components.

\vskip0.3cm

Naturally, at this point there are no guarantees that this procedure
performs well, let alone that it is close to optimal. That requires some
assumptions on the class $\F$, the hierarchy $(\F_\ell)_{\ell=1}^K$
and the regularization function $\Psi$. Moreover, the parameters that
each phase requires as inputs have to be specified for the procedure
to make any sense. All these issues are explored in what follows.
Before we dive into technicalities, and to give the reader a feeling
of how the regularized tournament looks like in a familiar situation,
let us describe the tournament version of the {\sc lasso}.

\subsection{The tournament {\sc lasso}} \label{sec-lasso-intro}
Consider the following standard setup: Let $X$ be an isotropic random vector in $\R^d$
(that is, for every $t \in \R^d$, $\E \inr{t,X}^2 = \|t\|_2^2$, where $\| \cdot \|_2$ denotes the Euclidean norm in $\R^d$).
Let $Y$ be the unknown target random variable and set $t_0$ to be the minimizer in $\R^d$ of the risk functional $t \to \E (\inr{X,t}-Y)^2$. For the sake of simplicity we assume that $Y=\inr{t_0,X} +W$ for $W$ that is mean-zero, square integrable and independent of $X$.

In sparse recovery problems one believes that $t_0$ is supported on at most $s$ coordinates with respect to the standard basis in $\R^d$---or at least it is well approximated by an $s$-sparse vector---, but one does not know that for certain. The {\sc lasso} procedure (introduced in \cite{Tib96}) selects $\wh{t}\in \R^d$ that minimizes the regularized empirical squared-loss functional
$$
t \to \frac{1}{N}\sum_{i=1}^N (\inr{t,X_i}-Y_i)^2 + \lambda \|t\|_1
$$
for a well chosen regularization parameter $\lambda$, and $\|t\|_1 =
\sum_{i=1}^d |t_i|$ is the $\ell_1$-norm of $t$.

The problem with the {\sc lasso} is that when either the class members $\inr{t,X}$ or the target $Y$ are heavy-tailed, the tradeoff between the accuracy with which the {\sc lasso} performs and the confidence with which that accuracy is attained is far from optimal. That suboptimal tradeoff is what the tournament {\sc lasso} aims to remedy.

For the time being we assume that for every $t \in \R^d$ and for every $p \leq c\log d$,  $\|\inr{X,t}\|_{L_p} \leq L \sqrt{p} \|\inr{X,t}\|_{L_2} = L\sqrt{p}\|t\|_2$ (where the last inequality holds because $X$ is isotropic). The random variable $W$ is assumed to be square-integrable.

The natural class of functions in this context is $\F = \{\inr{t,\cdot} : t \in \R^d\}$ and the hierarchy is given in terms of the sparsity of $t$: we set,
for $\ell=1,\ldots,K=\lceil\log_2 d\rceil+1$,
$$
\F_\ell = \{\inr{t,\cdot} : t \ {\rm is} \ d/2^{\ell-1} \ {\rm sparse}\}~.
$$
Let us go through the four phases of the tournament {\sc lasso}.
\subsubsection*{The first phase}
For the definition of the `referee' ${\cal DO}_\ell$ one requires to specify three parameters: $\rho_\ell$, $r_{\ell,1}$ and $n_1$. First, let $n_1 \sim_L N$, meaning that the cardinality of each block is a constant, depending only on $L$ and which we denote by $m_1$. Thus, $n_1=N/m_1$.

In what follows, $C(L),C_0(L),\ldots,C_4(L)$ denote appropriately
chosen constants whose value depends only on $L$. (The precise form
may be extracted from the analysis but it is of secondary importance
for our purpose.)

Next, set $s = d/2^{\ell-1}$. If $s \geq C(L) N/\log(ed/N)$ then the class $\F_\ell$ for the corresponding $\ell$ is too large and no useful statistical information can be derived from the sample. In that case, we set $r_{\ell,1},\rho_\ell=\infty$, and in particular, ${\cal DO}_\ell \equiv 0$. Otherwise, when $s <C_0(L) N/\log(ed/N)$, set
\begin{equation} \label{eq:LASSO-r-l-1}
r_{\ell} =C_1(L) \|W\|_{L_2} \sqrt{\frac{s}{d}\log\left(\frac{ed}{s}\right)} =C_2(L) \|W\|_{L_2} \sqrt{\frac{d \ell}{2^\ell N}},
\end{equation}
and
\begin{equation} \label{eq:LASSO-rho-l}
\rho_\ell = C_3(L) \|W\|_{L_2} \frac{s}{\sqrt{N}} \sqrt{\log\left(\frac{ed}{s}\right)} =C_4(L) \|W\|_{L_2} \frac{d \sqrt{\ell}}{2^\ell \sqrt{N}}.
\end{equation}
Consider $\wh{r}_\ell \geq r_\ell$ and set $r_{\ell,1} =\wh{r}_\ell$.

The first phase is performed as follows: given $(X_i)_{i=1}^N$, then for any $t_1,t_2$ supported on at most $d/2^{\ell-1}$ coordinates, let
\begin{description}
\item{$\bullet$} For $1 \leq j \leq n_1$ (which is proportional to $N$), set $v_j=\frac{1}{m_1} \sum_{i \in I_j} |\inr{t_1-t_2,X_i}|$ and define $d_\ell(t_1,t_2)$ to be a median of $\{v_j : 1 \leq j \leq n_1\}$.
\item{$\bullet$} Set ${\cal DO}_\ell(t_1,t_2) =1$ if either
$$
\|t_1-t_2\|_1 \geq \rho_\ell \ \ \ {\rm or} \ \ \
\|t_1-t_2\|_1 < \rho_\ell \ \  {\rm and} \ \ d_{\ell}(t_1,t_2) \geq r_{\ell,1}.
$$
\end{description}

\begin{Remark}
Since $X$ is isotropic, the $L_2(\mu)$ distance coincides with the
$\ell_2$ distance. Thus, ${\cal DO}_\ell$ can be chosen as a
deterministic function: ${\cal DO}_\ell(t_1,t_2)=1$ if $\|t_1-t_2\|_1
\geq \rho_\ell$ or if $\|t_1-t_2\|_1 < \rho_\ell$ and $\|t_1-t_2\|_2
\geq r_{\ell,1}$. However, when $X$ is not isotropic, the first phase
of the procedure is truly required, which is why we ignored the
simpler option that is available in the isotropic case.
\end{Remark}

\subsubsection*{The $\ell$-elimination phase}
The elimination phase in the tournament {\sc lasso} requires two parameters, $\lambda_\ell$ and $n_{\ell,2}$ as inputs, as well as the outcomes of ${\cal DO}_\ell$ obtained in the first phase. Recalling the choices of $r_{\ell,1}$ from \eqref{eq:LASSO-r-l-1} and $\rho_\ell$ from \eqref{eq:LASSO-rho-l}, we set
\begin{equation} \label{eq:LASSO-l-elimination}
r = c_1(L) r_{\ell,1} \ \  \lambda_\ell = c_2(L) \frac{r^2}{\rho_\ell^2}, \ \ {\rm and} \ \ n_{\ell,2}=c_3(L) N \min\left\{\frac{r^2}{\|W\|_{L_2}^2},1\right\},
\end{equation}
where $c_1,c_2$ and $c_3$ are constants that depend only on $L$. Let $(I_j)$ be the natural partition on $\{1,\ldots,N\}$ to $n_{\ell,2}$ disjoint blocks of equal cardinality, denoted by $m_{\ell,2}$. It follows that given the sample $(X_i,Y_i)_{i=N+1}^{2N}$ and $t_1,t_2 \in \F_\ell$, we set $t_1 \succ t_2$ if
$$
\frac{1}{m_{\ell,2}} \sum_{i \in I_j} \left(\inr{X_i,t_2}-Y_i)^2 - (\inr{X_i,t_1}-Y_i)^2\right) + \lambda_\ell \left(\|t_2\|-\|t_1\|\right) > 0
$$
for a majority of the blocks $I_j$. Therefore, $\cH^\prime_\ell$ consists of all  $t_1\in \R^d$ that satisfy $t_1 \succ t_2$ for any $t_2 \in \F_\ell$ such that ${\cal DO}_\ell(t_1,t_2)=1$.

\subsubsection*{The $\ell$-champions league phase}
The $\ell$-champions league phase receives as input the set $\cH_\ell^\prime$ produced in the second phase. Let $n_{\ell,2}$ and $r$ be as in \eqref{eq:LASSO-l-elimination} and set $r_{\ell,3} = c(L)r$. Given $(X_i,Y_i)_{i=2N+1}^{3N}$ and $t_1,t_2 \in \cH_\ell^\prime$ then $t_1 \gg t_2$ if
$$
\frac{2}{m_{\ell,2}} \sum_{i \in I_j} \inr{X_i,t_2-t_1} \cdot (\inr{X_i,t_1}-Y_i) \geq -r^2_{\ell,3},
$$
and $\cH_\ell = \{t_1 \in \cH_\ell^\prime \ : \ t_1 \gg t_2 \ {\rm for \ every \ } t_2 \in \cH_\ell^\prime \}$.

\subsubsection*{The fourth phase}
Given the classes $(\cH_\ell)_{\ell=1}^K$, $\ell_1$ is the largest integer $\ell$ such that $\bigcap_{j \leq \ell} \cH_j \not = \emptyset$. We select $\wh{t}$ to be any element of $\bigcap_{j \leq \ell_1} \cH_j$.

\vskip0.4cm The following theorem, proved in Section
\ref{sec:lassoproof}, summarizes the performance of the tournament {\sc lasso}.

\begin{Theorem} \label{thm:LASSO-intro}
  For $L \geq 1$ there are constants $c_0,\ldots,c_5$ that depend only on $L$ such that the following holds. Let $X$ be an isotropic random vector in $\R^d$. Let $Y=\inr{t_0,X}+W$ where $t_0\in \R^d$
and $W$ is mean-zero, square-integrable, and independent of $X$. Assume that for every $t \in \R^d$ and any $1 \leq p \leq c\log d$,
$\|\inr{t,X}\|_{L_p} \leq L \sqrt{p}\|\inr{t,X}\|_{L_2}$. Assume further that there is $v \in \R^d$ that is $s$-sparse such that
$$
\|t_0-v\|_1 \leq c_1(L)\|W\|_{L_2} \cdot s \sqrt{\frac{\log(ed/s)}{N}}~.
$$
If $N \geq c_2(L)s \log(ed/s)$, and
$$
\wh{r} \geq c_3(L) \|W\|_{L_2} \sqrt{\frac{s}{N} \log\left(\frac{ed}{s}\right)}
$$
then with probability at least
\begin{equation} \label{eq:LASSO-intro-conf}
1-2\exp\left(-c_4(L)N \min\left\{1,\left(\frac{\wh{r}}{\|W\|_{L_2}}\right)^2\right\}\right)~,
\end{equation}
we have
\begin{align} \label{eq:LASSO-intro-est}
& \|\wh{t}-t_0\|_2 \leq c_5(L)\wh{r}, \ \ \ \ \|\wh{t}-t_0\|_1 \leq c_5(L)\|W\|_{L_2} s \sqrt{\frac{\log(ed/s)}{N}}~,   \nonumber
\\
{\rm and} \ \ & \E(\inr{\wh{t},X}-Y)^2 \leq  \E(\inr{t_0,X}-Y)^2 +c_5(L) \wh{r}^2.
\end{align}
\end{Theorem}

Theorem \ref{thm:LASSO-intro} shows that the tournament {\sc
  lasso} attains the optimal accuracy/confidence tradeoff even though
$f^*(X)-Y$ can be heavy-tailed. In fact, the estimate is what one
would expect in the most friendly of scenarios: if $X$ were a
sub-Gaussian random vector (i.e., linear forms exhibiting a
$\psi_2-L_2$ norm equivalence with constant $L$ rather than an
$L_p-L_2$ moment equivalence going only up to $p \sim \log d$), and
$f^*(X)-Y$ were a Gaussian random variable, independent of $X$. The
{\sc lasso} does not come close to such an accuracy/confidence
tradeoff under the weak moment assumption of Theorem
\ref{thm:LASSO-intro}.

Note that $r =c_3 \|W\|_{L_2} \sqrt{\frac{s}{N}
  \log\left(\frac{ed}{s}\right)}$ is the best accuracy parameter one
can hope for even if the learner knows that $t_0$ is $s$-sparse and
$X$ and $W$ are Gaussian. The difference between the performance of
the tournament {\sc lasso} and the standard {\sc lasso} can be seen in
the confidence with which this accuracy is attained. In the situation
described in Theorem \ref{thm:LASSO-intro}, the standard {\sc lasso}
performs with that accuracy parameter only with constant confidence,
because all that we assume on $W$ is that it is square-integrable. In
contrast, the tournament {\sc lasso} attains the accuracy
\eqref{eq:LASSO-intro-est} with the optimal exponential probability
estimate \eqref{eq:LASSO-intro-conf}.

The tournament {\sc lasso} does not require prior information on the
degree of sparsity of $t_0$ to be carried out, but its success does
depend on having a large-enough sample and on that $\wh{r}$ is in the
right range. The former is a constraint that any recovery procedure
faces while the latter is easily achieved by running the procedure for
$\wh{r}_j=r_0/2^j$ for a large initial value of $r_0$, followed by a
standard validation argument at each step.

The one item that does require extra attention is that an upper estimate on $\|W||_{L_2}$ is used in the choice of parameters of the tournament {\sc lasso}. At times one is simply given that information; this is often the case in signal processing problems, where the nature of the `noise' is known to the learner. If not, one may use the data-dependent procedure from \cite{Men-tour3} which holds for more general noise models: it leads to upper and lower estimates on $\|f^*(X)-Y\|_{L_2}$ that are sharp up to absolute multiplicative constants and under minimal assumptions.

\vskip0.3cm

Of course, proving Theorem \ref{thm:LASSO-intro} requires some work,
and the choice of parameters used in the first three phases has to be
clarified. We explain the choice in the general case in the next two
sections and return to the example of the tournament {\sc lasso} in
Section \ref{sec:examples}.

\section{The main result}
\label{sec:main}

In the general setup we study, we merely assume a rather weak fourth-moment assumption.
More precisely, we work under the following conditions.

\begin{Assumption} \label{ass:moment}
Let $\F \subset L_2(\mu)$ be a locally compact, convex class of functions. Let $Y \in L_2$ and
assume that, for some constant $L>0$,
\begin{description}
\item{$\bullet$} for every $f,h \in \F$, $\|f-h\|_{L_4} \leq L \|f-h\|_{L_2}$;
\item{$\bullet$} $\|f^*-Y\|_{L_4} \leq \sigma_4$ for a known value $\sigma_4$.
\end{description}
\end{Assumption}

\begin{Remark}
The condition that $\|f^*-Y\|_{L_4} \leq \sigma_4$ may easily be
replaced by a combination of two assumptions: that for every $f \in
\F$, $\|f-Y\|_{L_4} \leq L\|f-Y\|_{L_2}$; and that $\|f^*-Y\|_{L_2}
\leq \sigma$ for some known constant $\sigma >0$. Also, in the case of
independent additive noise, that is, when $Y=f_0(X)+W$ where $f_0 \in
\F$ and $W$ that is mean-zero, square-integrable and independent of
$X$, the assumption that $\|f^*-Y\|_{L_4} \leq \sigma_4$ may be
replaced by the weaker one, that $\|W\|_{L_2} \leq \sigma$ for a known
constant $\sigma$.

The necessary modifications to the proofs are straightforward and we
do not explore this observation further. Also, as noted previously, we
refer the reader to \cite{Men-tour3} for a data-dependent procedure of
estimating $\|f^*-Y\|_{L_2}$ which may be easily modified to an
estimate on $\|f^*-Y\|_{L_4}$. Since that is not the main focus of
this paper we do not pursue it further and instead assume that the
learner has access to $\sigma_4$ or to $\sigma$.
\end{Remark}

\subsection{Complexity parameters of a class} \label{sec:complexity}
The choice of parameters $r_{\ell,1}$, $r_{\ell,3}$ and $\rho_\ell$ depends on a certain notion of ``complexity" of the underlying class. More accurately, there are geometric parameters that measure the metric entropy of ``localizations" of the class as well as the oscillation of various natural random processes indexed by those localizations.

The complexity is measured in terms of four parameters, depending both on the class $\F$ and the distribution of $(X,Y)$. The four play an essential role in describing the optimal performance of learning procedures and for detailed discussion on the meaning
we refer to Mendelson \cite{Men15,Men16} and Lugosi and Mendelson \cite{LuMe16}.

Before we define the four parameters we need some notation. Denote the unit ball in $L_2(\mu)$ by $D=\{f: \|f\|_{L_2}\le 1\}$ and let $S
=\{f: \|f\|_{L_2}= 1\}$ be the unit sphere. For
$h \in L_2(\mu)$ and $r>0$, we write $D_h(r) = \{f : \|f-h\|_{L_2} \leq
r\}$.
In a similar fashion for the norm $\Psi$ used as a regularization function, let ${\cal B} = \{f : \Psi(f) \leq 1\}$, set $\rho {\cal B} = \{f : \Psi(f) \leq \rho\}$ and  ${\cal B}_h(\rho) = \{f: \Psi(f-h) \leq \rho\}$.

In what follows we make two important modifications to the definitions of the complexity parameters used in \cite{Men15,Men16,LuMe16}. First, just like in the above-mentioned articles, we are interested in ``localized'' classes. However, because regularized procedures are affected by two norms, $\Psi$ and $L_2(\mu)$, the localization has to be with respect to both of them. Therefore, the ``localization'' of $\F$, centred in $h$ and of radii $\rho,r>0$ is defined by
$$
\F_{h,\rho,r} = (\F-h) \cap (\rho {\cal B} \cap rD)=\{f-h : f \in \F, \ \Psi(f-h) \leq \rho, \ \|f-h\|_{L_2} \leq r\}~.
$$

The second minor modification is that each complexity parameter is associated with the `worse case' centre $h \in \F^\prime$ for some fixed $\F^\prime \subset \F$, and not necessarily with the whole of $\F$.

\vskip0.4cm
Two of the four parameters are defined using the notion of \emph{packing numbers}.
\begin{Definition} \label{def:pack}
Given a set $H \subset L_2(\mu)$ and $\eps>0$, denote the
$\eps$-packing number of $H$ by
${\cal M}(H,\eps D)$. In other words, ${\cal M}(H,\eps D)$
is the maximal cardinality of a subset $\{h_1,\ldots,h_m\} \subset H$, for which $\|h_i-h_j\|_{L_2} \geq \eps$ for every $i \not = j$.
\end{Definition}
The first relevant parameter $\lambda_{\Q}$ is defined as
follows with appropriate numerical constants $\kappa$ and $\eta$:

\begin{Definition} \label{def-lambda-Q}
Fix $\rho>0$ and $h \in \F$. For $\kappa,\eta>0$, set
\begin{equation} \label{eq:covering-fixed-point-1}
\lambda_{\Q}(\kappa,\eta,h,\rho) =  \inf\{ r: \log {\cal M}(\F_{h,\rho,r},\eta r D) \leq \kappa^2 N\}~.
\end{equation}
For $\F^\prime \subset \F$ let
\[
\lambda_{\Q}(\kappa,\eta,\rho) = \sup_{h \in \F^\prime}
\lambda_{\Q}(\kappa,\eta,h,\rho)~.
\]
\end{Definition}
While $\kappa$ and $\eta$ are adjustable parameters, we are mainly interested in the behaviour of $\lambda_{\Q}$ as a function of $\rho$. The way one selects $\rho$ is clarified later.

The reason behind the choice of $\lambda_{\Q}$ comes from
high-dimensional geometry: if the class $\F_{h,\rho,r} \subset
L_2(\mu)$ is `less complex' than an $N$-dimensional Euclidean ball of
radius $\sim r$, then it can be covered by $\exp(cN)$ balls of radius
$\sim r$. Thus, $\lambda_{\Q}$ captures the smallest radius for which
there is still a chance that $\F_{h,\rho,r}$ resembles in some sense a
subset of a ball of radius $\sim r$ in $\R^N$. The hope is that for
such a choice of radius $r$, at least part of the metric structure of
$\F_{h,\rho,r}$ is reflected in a typical random set $\{
(u(X_i))_{i=1}^N : u \in \F_{h,\rho,r}\}$. In more statistical terms,
the hope is that a sample of cardinality $N$ provides the learner
with enough information to ``separate" class members that are far
enough, and the distance at which a sample of cardinality $N$ no
longer suffices is $\lambda_{\Q}$.

The next parameter, denoted by $\lambda_{\M}$, is also defined in terms
of the packing numbers of the localization $\F_{h,\rho,r}$, though at a different scaling than $\lambda_{\Q}$.
\begin{Definition} \label{def-lambda-M}
Fix $h \in \F$ and $\rho>0$. Let $\kappa>0$, $0<\eta<1$, and define
\begin{equation} \label{eq:covering-fixed-point-2}
\lambda_{\M}(\kappa,\eta,h,\rho) =  \inf\{ r: \log {\cal M}(\F_{h,\rho,r},\eta r D) \leq \kappa^2 N r^2\}~.
\end{equation}
Also, for $\F^\prime \subset \F$ let
\[
\lambda_{\M}(\kappa,\eta,\rho) = \sup_{h \in \F^\prime}
\lambda_{\M}(\kappa,\eta,h,\rho)~.
\]
\end{Definition}

For the remaining two complexity parameters, let $(\eps_i)_{i=1}^N$ be
independent, symmetric $\{-1,1\}$-valued random variables that are
independent of $(X_i,Y_i)_{i=1}^N$.

\begin{Definition} \label{def:fixed-emp}
Fix $h \in \F$ and $\rho>0$. For $\kappa>0$ let
\begin{equation} \label{eq:emp-fixed-point}
r_{E}(\kappa,h,\rho) = \inf\left\{ r: \E \sup_{u \in \F_{h,\rho,r}} \left| \frac{1}{\sqrt{N}} \sum_{i=1}^N \eps_i u(X_i) \right| \leq \kappa \sqrt{N} r\right\},
\end{equation}
and for $\F^\prime \subset \F$ set $r_{E}(\kappa,\rho) = \sup_{h \in \F^\prime} r_{E}(\kappa,h,\rho)$.
\end{Definition}
The idea behind $r_E$ is similar to the one behind $\lambda_{\Q}$. $r_E(\kappa,h,\rho)$ is the smallest radius for which the expected supremum of the process $u \to \frac{1}{\sqrt{N}} \sum_{i=1}^N \eps_i u(X_i)$ indexed by $\F_{h,\rho,r}$ exceeds the expectation of the supremum of the Bernoulli process in $\R^N$, $v \to \sum_{i=1}^N \eps_i v_i$, indexed by a Euclidean ball of radius $r$. Thus, the point is to identify when a typical random set $\{ N^{-1/2}(u(X_i))_{i=1}^N : u \in \F_{h,\rho,r}\}$ is richer than Euclidean ball of radius $\sim r$ in $\R^N$, where `richer' is measured in terms of the Bernoulli mean-width of a set rather than by its metric entropy.

\begin{Definition} \label{def:fixed-multi}
Fix $h \in \F$ and $\rho>0$. For $\kappa>0$, set $\ol{r}_{\M}(\kappa,h,\rho)$ to be
\begin{equation}
\ol{r}_{\M}(\kappa,h,\rho) = \inf\left\{ r: \E \sup_{u \in \F_{h,\rho,r}} \left| \frac{1}{\sqrt{N}} \sum_{i=1}^N\eps_i u(X_i) \cdot (h(X_i)-Y_i) \right| \leq \kappa \sqrt{N} r^2\right\}.
\end{equation}
For $\sigma>0$ put $\F_Y^{(\sigma)}=\{ f \in \F^\prime : \|f(X)-Y\|_{L_2} \leq \sigma\}$ and let
${\wt r}_{\M}(\kappa,\sigma,\rho) = \sup_{h \in \F_Y^{(\sigma)}} {\ol r}_{\M}(\kappa,h,\rho)$.
\end{Definition}

Finally, suppose that the distribution of $(X,Y)$ is such that
$\|Y-f^*(X)\|_{L_4} \leq \sigma_4$  for a known constant $\sigma>0$. The ``complexity'' of $\F$ relative to centres in $\F^\prime$ and radius $\rho$ is
\begin{equation} \label{eq:r*}
r^*(\F,\F^\prime,\rho) =\max\{\lambda_{\Q}(c_1,c_2,\rho),\lambda_{\M}(c_1/\sigma_4,c_2,\rho),r_E(c_1,\rho),\wt{r}_{\M}(c_1,\sigma_4,\rho) \}~.
\end{equation}
Here $c_1,c_2$ are appropriate positive numerical
constants. (``Appropriate'' means that $r^*(\F,\F^\prime,\rho)$ satisfies
Propositions \ref{thm:distance-functional}, \ref{thm:elimination} and \ref{thm:winners} below).
The existence of such constants is proved in \cite{LuMe16} when $\F^\prime = \F$, i.e., when any function in $\F$ is a `legal choice' of a centre, and under Assumption \ref{ass:moment}. In that case, the constants depend only on the value of $L$.

When $\F$ and $\F^\prime$ are clear from the context, we simply write $r^*(\rho)$ for $r^*(\F,\F^\prime,\rho)$.

\subsection{Properties of the hierarchy}
After the complexity parameters are set in place, let us identify the conditions on the hierarchy $(\F_\ell)_{\ell=1}^K$ that are needed in the analysis of the regularized tournament.

Recall that $\F$ is a (convex) subset of a normed space $(E,\Psi)$; $E$ is also a subspace of $L_2(\mu)$, though the norms $\Psi$ and $\| \cdot \|_{L_2(\mu)}$ may have nothing to do with each other.
Let $B_{\Psi^*}$ and $S_{\Psi^*}$ denote the unit ball and unit sphere in the dual space to $(E,\Psi)$, respectively. Therefore, $B_{\Psi^*}$ consists of all the linear functionals $z \in E^*$ for which $\sup_{\{x \in E : \Psi(x)=1\}} |z(x)| \leq 1$. A linear functional $z^* \in S_{\Psi^*}$ is a norming functional for $f \in E$ if $z^*(f)=\Psi(f)$.

\begin{Definition} \label{def:Delta}
Let $\Gamma_{f}(\rho) \subset S_{\Psi^*}$ be the collection of functionals that are norming for some $v \in {\cal B}_f(\rho/20)$. Set
\begin{equation*}
\Delta_\ell(\rho,r) = \inf_{f \in \F_\ell} \inf_h \sup_{z \in \Gamma_{f}(\rho)} z(h-f)~,
\end{equation*}
where the inner infimum is taken in the set
\begin{equation} \label{eq:condition-on-set}
\{ h \in \F: \Psi(h-f)=\rho \ {\rm and} \ \|h-f\|_{L_2} \leq r\}~.
\end{equation}
\end{Definition}


Let us examine $\Delta_\ell(\rho_\ell,r_\ell)$ and explain its meaning for some fixed values $\rho,r>0$. Note that $\Delta_\ell(\rho,r) \leq \rho$. Indeed, $\Gamma_{f}(\rho) \subset S_{\Psi^*}$ and if $z \in S_{\Psi^*}$ and $\Psi(h-f) \leq \rho$ then
$$
|z(h-f)| \leq \Psi^*(z) \cdot \Psi(h-f) \leq \rho~.
$$
The interesting situation is when one can ensure a reverse inequality, that is, that $\Delta_\ell(\rho,r)$ is proportional to $\rho$, say $\Delta_\ell(\rho,r) \geq (4/5)\rho$. Such a lower estimate on $\Delta_\ell$ implies the following. Let $f \in \F_\ell$ and $h \in \F$ for which
$\Psi(h-f)=\rho$ and $\|f-h\|_{L_2} \leq r$. It follows that there is some $z \in S_{\Psi^*}$ and $v \in {\cal B}_f(\rho/20)$ such that $z$ is norming for $v$ and $z(h)-z(f) \geq \Delta_\ell(\rho,r)$. Therefore,
\begin{eqnarray*}
\Psi(h)-\Psi(f) & = & \Psi(h)-\Psi(v+(f-v)) \geq \Psi(h)-\Psi(v)-\Psi(f-v)
\\
& \geq & z(h)-z(v)-\Psi(f-v) \geq z(h)-z(f) - 2\Psi(f-v)
\\
& \geq & \Delta_\ell(\rho,r) -\rho/10 \geq 3\rho/5~.
\end{eqnarray*}
A lower bound on the term $\Psi(h)-\Psi(f)$ plays an essential role in the study of the elimination phase of the regularized tournament, when one has to compare
$$
\frac{1}{m_{\ell,2}}\sum_{i \in I_j} (f(X_i)-Y_i)^2 + \lambda_\ell \Psi(f) \ \ {\rm and} \ \ \frac{1}{m_{\ell,2}}\sum_{i \in I_j} (h(X_i)-Y_i)^2 + \lambda_\ell \Psi(h).
$$

Obviously, ensuring that $\Delta_{\ell}(\rho,r) \geq (4/5)\rho$ becomes simpler when the set $\Gamma_f(\rho)$ is large.  In the extreme case, when $\rho > 30\Psi(f)$, then ${\cal B}_f(\rho/20)$ contains a nontrivial $\Psi$-ball around $0$; thus,
$\Gamma_f(\rho)=S_{\Psi^*}$ and $\Delta_\ell(\rho,r)=\rho$. The other extreme is if $\rho$ is very small and one is left only with the functionals that are norming for $f$ itself. Intuitively, the right choice of $\rho_\ell$ is the smallest one for which, for $r_\ell=r^*(\F,\F_\ell,\rho_\ell)$, one has $\Delta_\ell(\rho_\ell,r_\ell) \geq 4\rho_\ell/5$.
\begin{framed}
\begin{Definition} \label{def:compatible}
The sequence $(\F_\ell,\rho_\ell)_{\ell=1}^K$ is compatible if
\begin{description}
\item{(1)}  $\F=\F_1 \supset \F_2 \supset \cdots \supset \F_K$ is a finite hierarchy;
\item{(2)} $(\rho_\ell)_{\ell=1}^K$ is decreasing and $r_\ell=r^*(\F,\F_\ell,\rho_\ell)$;
\item{(3)} for every $1 \leq \ell \leq K$, $\Delta_\ell(\rho_\ell,r_\ell) \geq 4\rho_\ell/5$.
\end{description}
We allow the choice of $\rho_\ell=r_\ell=\infty$ for $\ell=1,\ldots,\ell_0$. In such cases the compatibility condition is to be verified from $\ell_0+1$ onward.
\end{Definition}
\end{framed}

We are now ready to specify the parameters used in the definition of a regularized tournament.

\begin{framed}
\begin{description}
\item{$\bullet$} Let $\alpha$, $\beta$, $m_1$, $\theta_1$ and $\theta_2$ be well chosen constants that depend only on the norm equivalence constant $L$ from Assumption \ref{ass:moment}, and assume that one has access to the value $\sigma_4$ from that assumption.
\item{$\bullet$} Assume further that $(\F_\ell,\rho_\ell)_{\ell=1}^K$ is a compatible sequence, let $r_\ell$ be as in Definition \ref{def:compatible} and set $\wh{r}_\ell > r_\ell$.
\end{description}
We set the following choice of parameters:
\begin{description}
\item{\underline{\rm First phase:}} Let $\rho_\ell$ as above, and set
$$
n_1 = \frac{N}{m_1} \ \ \ \ {\rm and} \ \ \ \ r_{\ell,1} = \beta \wh{r}_\ell.
$$
\item{\underline{\rm $\ell$-elimination phase:}} Set
$$
n_{\ell,2} = \theta_1 N \min\left\{1,\frac{\wh{r}_\ell^2}{\sigma_4^2}\right\} \ \ {\rm and} \ \ \lambda_\ell=\theta_2 \frac{\wh{r}_\ell^2}{\rho_\ell}.
$$
\item{\underline{\rm $\ell$-champions league:}} Set $r_{\ell,3} = c(\beta/\alpha)\wh{r}_\ell$ for a suitable absolute constant $c$ ($c=1/\sqrt{5}$ would do).
\end{description}
\end{framed}

With these choices set in place, let us formulate the main result of this article.
\begin{Theorem} \label{thm:main}
For $L \geq 1$ there exist constants $c_0,\ldots,c_2$ that depend only on $L$ and for which the following holds. Let $(\F,X,Y)$ satisfy Assumption \ref{ass:moment} and set $(\F_\ell)_{\ell=1}^K$ be a compatible sequence of the class $\F$. If  $\ell^*$ is the largest index $\ell$ such that $f^* \in \F_\ell$, then  with probability at least
$$
1-2 \sum_{\ell=1}^{\ell^*} \exp\left(-c_0(L)N \min\left\{1,\frac{\wh{r}_{\ell}^2}{\sigma_4^2}\right\}\right),
$$
we have
$$
\Psi(\wh{h}-f^*) \leq \rho_{\ell^*}, \ \ \|\wh{h}-f^*\|_{L_2} \leq c_1(L)\wh{r}_{\ell^*}, \ \ {\rm and} \ \ R(\wh{h})-R(f^*) \leq c_2(L) \wh{r}_{\ell^*}^2.
$$
\end{Theorem}

\subsection{Discussion}


The main message of Theorem \ref{thm:main} is that the regularized
tournament procedure achieves essentially the best performance that
can be expected even under strong assumptions of sub-Gaussian distributions.
The regularized procedure yields (almost) the optimal
accuracy-confidence tradeoff for any accuracy parameter $r \geq
r_{\ell^*}$: it behaves as if it ``knew'' the location of $f^*$ in the
hierarchy without actually knowing it. Indeed, the accuracy/confidence
tradeoff established in Theorem \ref{thm:main} is essentially the best
possible for any learning procedure taking values only in
$\F_{\ell^*}$.

We emphasize that Theorem \ref{thm:main} is quite general though
finding the adequate parameters of the procedure requires additional work.
For the ``tournament'' version of {\sc lasso} and {\sc slope} we
work out the details in Section \ref{sec:examples} under certain assumptions
(such as isotropic design vector $X$ and approximately sparse
linear regression function) for illustration. Some of these assumption may
be weakened but we prefer to keep the presentation as simple as possible.

\noindent
{\bf Related work.}
The sensitivity of empirical risk minimization (or least squares regression)
to heavy-tailed distributions has been pointed out and several proposals
of robust regression function estimates have been made that avoid
this sensitivity. We refer to
Audibert and Catoni \cite{AuCa11},
Hsu and Sabato \cite{HsuSa13},
Lerasle and Oliveira \cite{LeOl12},
Minsker \cite{Min15},
Brownlees, Joly, and Lugosi \cite{BrJoLu15},
Lugosi and Mendelson \cite{LuMe16}
for a sample of the literature. This paper mostly builds upon
the methodology of median-of-means tournaments, developed in \cite{LuMe16},
(see also Lugosi and Mendelson \cite{LuMe17}).
Here we extend this methodology to the analysis of regularized
robust risk minimization similarly to how the paper of
Lecu\'e and Mendelson \cite{LeMe16} analyzes standard regularized
risk minimization. The analysis of the {\sc lasso} and {\sc slope}
procedures of \cite{LeMe16} was extended and generalized by
Bellec, Lecu{\'e}, and Tsybakov \cite{BeLeTs16}.
In an independent parallel work to ours, and building on the arguments developed in \cite{LuMe16}, Lecu{\'e} and Lerasle \cite{LeLe17}
point out a connection of median-of-means tournaments to
Le Cam's estimators, develop a version of {\sc lasso}--the
so-called {\sc mom-lasso}--and prove a performance bound quite similar
to Theorem \ref{thm:LASSO-intro}.

\section{Analyzing the four phases}
\label{sec:proof}

\subsection*{The first phase - the $\ell$-distance oracle}
At each stage $\ell=1,\ldots,K$ of the procedure, one initially uses a modification of the distance oracle from \cite{LuMe16}.

The $\ell$-distance oracle is a data-dependent procedure that provides information on the distances between functions. It is used for any pair $f,h \in \F$, and aims at determining if $\Psi(f-h) \geq \rho_\ell$, or, if  $\Psi(f-h) \leq \rho_\ell$, whether $\|f-h\|_{L_2}\geq r_{\ell,1}$. Note that $\Psi$ is a known norm and therefore, $\Psi(f-h)$ is known for any pair $f,h\in \F$ but $\|f-h\|_{L_2}$ needs to be (crudely) estimated.

Recall that we work under Assumption \ref{ass:moment}, and that ${\cal DO}_\ell$ is defined as follows: we split $\{1,\ldots,N\}$ to $n=n_1$ disjoint blocks $(I_j)_{j=1}^n$, each one of cardinality
$m=N/n$. For a sample $\C_1=(X_i)_{i=1}^N$ and functions $f$ and $h$, let $w=(|f(X_i)-h(X_i)|)_{i=1}^N$ and set
$$
\Phi_{\C_1}(f,h) = {\rm Med}_m(w)~,
$$
where ${\rm Med}_m(w)$ is a median of the $n$ values $\frac{1}{m}\sum_{i \in I_j} |(f-h)(X_i)|$.

The behaviour of $\Phi$ described below has been established in Mendelson \cite{Men16a} (see also Lugosi and Mendelson \cite{LuMe16}):
\begin{Proposition} \label{thm:distance-functional}
Let $\F$ satisfy Assumption \ref{ass:moment}. There exist constants $m$, $0<\alpha < 1 <
\beta$, and $\kappa,\eta$ and $c$, all of them depending only on $L$, for which the following holds.
\begin{description}
\item{$\bullet$} For $1 \leq \ell \leq K$ and $\rho>0$, let $r>r^*(\F,\F_\ell,\rho)$, where $r^*$ is defined relative to the constants $\kappa$ and $\eta$).
\item{$\bullet$} Let $n=N/m$ and fix $f \in \F_\ell$.
\end{description}

Then, with probability at least $1-2\exp(-cN)$, for any $h \in \F$ that satisfies $\Psi(f,h) \leq \rho$,
\begin{description}
\item{$(1)$} if $\Phi_{\C_1}(f,h) \geq \beta r$ then
$$
\beta^{-1} \Phi_{\C_1}(f,h) \leq \|f-h\|_{L_2} \leq \alpha^{-1} \Phi_{\C_1}(f,h);
$$
\item{$(2)$} if $\Phi_{\C_1}(f,h) < \beta r$ then $\|f-h\|_{L_2} \leq (\beta/\alpha)r$.
\end{description}
\end{Proposition}

The proof of Proposition \ref{thm:distance-functional} is a direct outcome of Proposition 3.2 from \cite{LuMe16}, applied to the set $\F \cap {\cal B}_{f}(\rho)$ for a fixed centre $f \in \F_\ell$ and we shall not present it here.

\vskip0.3cm

Based on Proposition \ref{thm:distance-functional}, our choice of parameters in the first phase is clear. Recall that for $\wh{r}_\ell > r_\ell=r^*(\F,\F_\ell,\rho_\ell)$ and $f_1,f_2 \in \F$ we set ${\cal DO}_\ell(f_1,f_2)=1$ if either $\Psi(f_1-f_2) > \rho_\ell$ or if $\Psi(f_1-f_2) \leq \rho_{\ell}$ and $\Phi_{\C_1}(f_1,f_2) \geq \beta \wh{r}_\ell$.

\begin{framed}
Thanks to Proposition \ref{thm:distance-functional}, it follows that for a fixed centre $f \in \F_\ell$ (which is selected as $f^*$ in what follows), with probability at least $1-2\exp(-cN)$ if $h \in \F$, and ${\cal DO}_\ell(f,h) =0 $ then $\Psi(h-f) \leq \rho_\ell$ and $\|f-h\|_{L_2} \leq \wh{r}_\ell$.
\end{framed}

\begin{Remark}
Although Proposition \ref{thm:distance-functional} is formulated for a designated single centre $f$, it is straightforward to extend it to \emph{any} centre in $\F$ and obtain a uniform distance oracle that holds for any pair $f,h \in \F$.
\end{Remark}

\subsection*{The second phase---$\ell$-elimination}
Fix $1 \leq \ell \leq K$ and let $f,h \in \F$.
Recall the definition of a regularized match between $f$ and $h$:
first, the
$\ell$-distance oracle defined above uses the first part of the sample
${\cal C}_1=(X_i,Y_i)_{i=1}^N$ to determine the value of ${\cal DO}_\ell(f,h)$.  If ${\cal DO}_\ell(f,h)=0$, the match is abandoned.

Each match that is allowed to take place by the $\ell$-distance oracle is played using the second part of the sample,
${\cal C}_2=(X_i,Y_i)_{i=N+1}^{2N}$. The sub-sample is partitioned to $n=n_{\ell,2}$ blocks $(I_j)_{j=1}^n$ of cardinality $m=N/n$ where recall that $n$ is chosen as
\begin{equation} \label{eq:n-l-2-body}
\theta_1(L) N\min\left\{1,\frac{\wh{r}_\ell^2}{\sigma_4^2}\right\}.
\end{equation}
for a well-chosen constant $\theta_1$ that depends only on the equivalence constant $L$ from Assumption \ref{ass:moment}. We also set
$$
\lambda_\ell=\theta_2(L) \frac{\wh{r}_\ell^2}{\rho_\ell},
$$
with the choices of both constants $\theta_1$ and $\theta_2$ specified below.

The key definition in the elimination stage is the choice of a winner in a `statistical match' between two functions in $\F_\ell$.
\begin{Definition}
The function $f$ defeats $h$ (denoted by $f \succ h)$ if
$$
\frac{1}{m}\sum_{i \in I_j} \left((h(X_i)-Y_i)^2 - (f(X_i)-Y_i)^2\right) + \lambda_\ell (\Psi(h)-\Psi(f)) > 0
$$
on a majority of the blocks $I_j$.

The set of winners $\cH^\prime_\ell$ of the $\ell$-th elimination round consists of all the functions in $\F_\ell$ that have not lost a single match against a function in $\F_\ell$.
\end{Definition}

The idea behind the elimination round is to `exclude' functions that are far from $f^*$ (without knowing the identity of $f^*$, of course). To that end, it suffices to show that, if $f^* \in \F_\ell$, then with high probability, $f^*$ wins all the matches it takes part in. Indeed, that implies that $\cH_\ell^\prime$ is nonempty, and that all the matches between $f^*$ and any $h \in \cH^\prime_\ell$ must have been abandoned; therefore, ${\cal DO}_\ell(f^*,h)=0$, that is,
\begin{equation} \label{eq:prelim-phase-outcome}
\Psi(f^*-h) \leq \rho_\ell \ \ {\rm and} \ \ \|f^*-h\|_{L_2} \leq (\beta/\alpha)\wh{r}_\ell~.
\end{equation}
The next theorem describes the outcome of the $\ell$-elimination phase. Its proof may be found in Section \ref{sec:proof-of-elimination}.

\begin{Proposition} \label{thm:elimination}
Using the notation above, if $f^* \in \F_\ell$ then, with probability at least
$$
1-2\exp\left(-c_0N \min\left\{1,\frac{\wh{r}_\ell^2}{\sigma_4^{2}}\right\}\right)~,
$$
$f^*$ wins all the matches is participates in. In particular, on that event, if $h \in \cH^\prime_\ell$, then \eqref{eq:prelim-phase-outcome} holds.
\end{Proposition}

\subsection{Proof of Proposition \ref{thm:elimination}---highlights} \label{sec:proof-of-elimination}

To explain why this elimination phase preforms well even when $\F$ is
very large, define, for each block $I_j$ ($j=1,\ldots,n$),
\begin{equation} \label{eq:reg-emp-excess}
B^\lambda_{h,f}(j)=\frac{1}{m}\sum_{i \in I_j} \left((h(X_i)-Y_i)^2 - (f(X_i)-Y_i)^2\right) + \lambda (\Psi(h)-\Psi(f))~.
\end{equation}
Note that the regularized empirical excess risk of $h$ on block $I_j$ is $B^\lambda_{h,f^*}(j)$.

Consider the $\ell$-th stage of the regularized tournament.  The assertion of Proposition \ref{thm:elimination} is that, if $f^* \in \F_\ell$, then it is a winner of all the elimination phase matches it participates in. Hence, Proposition \ref{thm:elimination} is proved once we ensure that for the right choice of $\lambda=\lambda_\ell$, with high probability, if $h \in \F$ and ${\cal DO}_\ell(f^*,h)=1$ then $B^\lambda_{h,f^*}(j)$ is positive for most of the blocks $I_j$.

To that end, observe that
\begin{eqnarray*}
\lefteqn{
\frac{1}{m}\sum_{i \in I_j} \left((h(X_i)-Y_i)^2 -
  (f^*(X_i)-Y_i)^2\right) }
\\
& = &
\frac{1}{m}\sum_{i \in I_j} (h-f^*)^2(X_i) + \frac{2}{m} \sum_{i \in I_j} (h-f^*)(X_i) \cdot (f^*(X_i)-Y_i)~,
\end{eqnarray*}
which is the natural decomposition of the empirical excess risk functional into its quadratic and multiplier components. Setting
$$
\Q_{h,f}(j)=\frac{1}{m}\sum_{i \in I_j} (h-f)^2(X_i) \ \ {\rm and} \ \ \M_{h,f}(j) = \frac{2}{m} \sum_{i \in I_j} (h-f)(X_i) \cdot (f(X_i)-Y_i)~,
$$
we have
\begin{eqnarray*}
B^\lambda_{h,f^*}(j) & = & \frac{1}{m}\sum_{i \in I_j} \left((h(X_i)-Y_i)^2 - (f^*(X_i)-Y_i)^2\right) + \lambda(\Psi(h)-\Psi(f^*))
\\
& = & \Q_{h,f^*}(j) + \M_{h,f^*}(j) + \lambda(\Psi(h)-\Psi(f^*))~.
\end{eqnarray*}

The first observation we require is a version of a deterministic
result from \cite[Theorem 3.2]{LeMe16a} (see the appendix for the proof).

\begin{Lemma} \label{lemma:basic-combining-loss-and-reg}
Let $f^* \in \F_\ell$ and $h \in \F$ for which either $\Psi(h-f^*) = \rho$, or $\Psi(h-f^*) < \rho$ and $\|h-f^*\|_{L_2} \geq r$. Assume that $\Delta_\ell(\rho,r) \geq 4\rho/5$, that $\lambda$ satisfies
\begin{equation} \label{eq:lambda}
\frac{C}{2} \cdot \frac{r^2}{\rho} \leq \lambda \leq \frac{3C}{4} \cdot \frac{r^2}{\rho}
\end{equation}
for some $C>0$. Assume further that
\begin{equation} \label{eq:cond1}
\M_{h,f^*}(j) \geq -(C/4)\max\left\{\|h-f^*\|_{L_2}^2,r^2\right\},
\end{equation}
and if also $\|h-f^*\|_{L_2} \geq r$ then
\begin{equation} \label{eq:cond2}
\Q_{h,f^*}(j) \geq C \|h-f^*\|_{L_2}^2.
\end{equation}
Then
$$
\Q_{h,f^*}(j) + \M_{h,f^*}(j) + \lambda(\Psi(h)-\Psi(f^*)) > 0~.
$$
Moreover, if $f \in \F$ such that $f=f^*+\alpha(h-f^*)$ for $\alpha>1$ then also
$$
\Q_{f,f^*}(j) + \M_{f,f^*}(j) + \lambda(\Psi(f)-\Psi(f^*)) > 0~.
$$
\end{Lemma}

Thanks to Lemma \ref{lemma:basic-combining-loss-and-reg} (which is proved in the appendix), all that is required to prove Proposition \ref{thm:elimination} is to verify that
with the requested probability, \eqref{eq:cond1} and \eqref{eq:cond2}
hold uniformly in $h$ and on a majority of the blocks $I_j$, provided
that $f^* \in \F_\ell$; thus $f^*$ defeats any function in $h \in \F$ that satisfies either $\Psi(h-f^*) \geq \rho$ or $\Psi(h-f^*) \leq \rho$ and $\|f^*-h\|_{L_2} \leq r$.

The fact that we have the required control over coordinate blocks is formulated in the following lemma. Its proof may be found in the appendix.

\begin{Lemma} \label{lemma:components-of-main}
There exists an absolute constant $c$ and a constant $C_1=C_1(L,\tau)$
for which the following holds. Let $f^* \in \F_\ell$. For $0<\tau<1$, with probability at least $1-2\exp(-c\tau^2 n)$, for every $h \in {\cal B}_{f^*}(\rho)$ that satisfies $\|f-f^*\|_{L_2} \geq \wh{r}_\ell$, we have
$$
\left|\left\{j : {\Q}_{f,f^*}(j) \geq C_1 \|f-f^*\|_{L_2}^2 \right\} \right| \geq \left(1-\tau\right)n
$$
and
$$
\left|\left\{j : {\M}_{f,f^*}(j) \leq -\frac{C_1}{4}\|f-f^*\|_{L_2}^2\right\} \right| \leq \tau n~.
$$
Moreover, for every $h \in {\cal B}_{f^*}(\rho)$ that satisfies $\|f-f^*\|_{L_2} \leq \wh{r}_\ell$, we have
$$
\left|\left\{j : {\M}_{f,f^*}(j) \leq -\frac{C_1}{4}r^2\right\} \right| \leq \tau n~.
$$
\end{Lemma}

It follows from Lemma \ref{lemma:components-of-main} that if $\tau <1/4$ then with probability at least $1-2\exp(-c\tau^2 n)$, for every $h$ as in Lemma \ref{lemma:basic-combining-loss-and-reg}, conditions \eqref{eq:cond1} and \eqref{eq:cond2} hold for $C=C_1$ and $r=\wh{r}_\ell$ on the majority of the blocks $I_j$. Hence, setting $\tau=1/10$, on an event with probability at least $1-2\exp(-cn)$, $f^*$ wins all the matches it participates in and that are allowed to take place by the $\ell$-distance oracle, as we require.

\subsection*{The third phase---$\ell$-champions league} \label{sec:champ}

Once Proposition \ref{thm:elimination} is established, we turn to the
third phase, aimed at selecting a set of ``winners". In order to do
that, we run the champions league tournament defined in Definition \ref{def:phase-3}, performed in each one of the sets
$\cH^\prime_\ell$. (This is the same procedure used in \cite{LuMe16}.)
The crucial point is that if $f^* \in \F_\ell$ then
$\cH_\ell^\prime$ satisfies the necessary conditions for a champions
league tournament: that $f^*\in \cH_\ell^\prime$ and all the functions
$h\in \cH_\ell^\prime$ have a mean-squared error at most $ \sim
\wh{r}_\ell$.

Recall that the $\ell$-champions league consists of matches that use
the third part of the sample ${\cal C}_3=(X_i,Y_i)_{i=2N+1}^{3N}$. Let
$(I_j)_{j=1}^n$ be the partition of $\{2N+1,\ldots,3N\}$ to $n$
blocks, for the same value of $n=n_{\ell,2}$ as in the
$\ell$-elimination phase, given in \eqref{eq:n-l-2-body}. Also, for
$\alpha$ and $\beta$ as in Proposition \ref{thm:distance-functional},
set $c=\beta/\alpha$ and for $f,h \in \cH_\ell^\prime$, let
$\Psi_{h,f} = (h (X)-f (X))(f(X)-Y)$.

Recall that the function $f$ wins its home match against $h$ (denoted by $f \gg h$) if
$$
\frac{2}{m} \sum_{i \in I_j} \Psi_{h,f}(X_i,Y_i) \geq -(2c \wh{r}_\ell)^2/10
$$
on more than $n/2$ of the blocks $I_j$, and the set of winners $\cH_\ell$ consists of all the ``champions" in $\cH_\ell^\prime$ that win all of their home matches.

\vskip0.3cm

The outcome of the $\ell$-champions league phase is as follows:
\begin{Proposition} \label{thm:winners}
Let $\cH_\ell^\prime$ as above. With probability at least
$$
1-2\exp\left(-c_0N \min\left\{1,\frac{\wh{r}_\ell^2}{\sigma_4^{2}}\right\}\right)
$$
with respect to $(X_i,Y_i)_{i=2N+1}^{3N}$, the set of winners $\cH_\ell$ contains $f^*$, and if $h \in \cH_\ell$ then
$$
R(h) - R(f^*) \leq 16c^2\wh{r}_\ell^2~.
$$
\end{Proposition}
Proposition \ref{thm:winners} is an immediate outcome of Proposition 3.8 from \cite{LuMe16} for $H=\cH_\ell^\prime$ and using the fact that $f^* \in \cH_\ell^\prime$ and that if $h \in \cH_\ell^\prime$ then $\|h-f^*\|_{L_2} \leq (\beta/\alpha)\wh{r}_\ell$.

Combining all these observations, Corollary \ref{cor:sum} describes the outcome of the first three phases in the regularized tournament procedure, for each member of the hierarchy.

\begin{Corollary} \label{cor:sum}
For $L \geq 1$ there exists a constant $c$ that depends only on $L$ for which the following holds. Using the above notation, if $f^* \in \F_\ell$, then with probability at least
$$
1-2\exp\left(-c_1(L)N \min\left\{1,\frac{\wh{r}_\ell^2}{\sigma_4^2}\right\}\right)
$$
with respect to $(X_i,Y_i)_{i=1}^{3N}$, the set of winners $\cH_\ell$ satisfies:
\begin{description}
\item{$\bullet$} $f^* \in \cH_\ell$, and
\item{$\bullet$} for any $h \in \cH_\ell$,
$$
\Psi(h-f^*) \leq \rho_\ell, \ \ \|h-f^*\|_{L_2} \leq (\beta/\alpha)\wh{r}_\ell, \ \ {\rm and} \ \ R(h)-R(f^*) \leq 16(\beta/\alpha)^2 \wh{r}_\ell^2~.
$$
\end{description}
\end{Corollary}

\subsection*{Selection of a final winner}

Once the first three phases are completed, their outcome is $K$ sets
$\cH_\ell$ consisting of stage winners. Of course, some of these sets
may be empty. However, for those indices $\ell$ for which $f^*\in
\F_\ell$, with high probability the set $\cH_\ell$ contains at least
$f^*$. In order to select the final ``winner'', let $\ell_1$ be the
largest integer $1 \leq \ell \leq K$ for which $\bigcap_{j \leq \ell}
\cH_j \not = \emptyset$ and that the procedure returns any $\wh{h} \in
\bigcap_{j \leq \ell_1} \cH_j$.

Clearly, on a high-probability event, $\ell_1 \geq \ell^*$.
On this event, the selected function $\wh{h}$ belongs to $\cH_{\ell^*}$. Moreover, recalling that $(\wh{r}_\ell)_{\ell=1}^K$ is decreasing,
$$
\Psi(\wh{h}-f^*) \leq \wh{r}_{\ell^*}, \ \ \|\wh{h}-f^*\|_{L_2} \leq (\beta/\alpha)\wh{r}_{\ell^*} \ \ {\rm and} \ \
R(\wh{h})-R(f^*) \leq c(\beta/\alpha)^2 \wh{r}_{\ell^*}^2~,
$$
which completes the proof of Theorem \ref{thm:main}.
\endproof

\section{Examples} \label{sec:examples}
In what follows we present two examples: a tournament version of {\sc lasso}, and also, a tournament version of another popular sparse recovery procedure---{\sc slope}.

Here we use a stronger moment assumption than Assumption \ref{ass:moment} because that allows us to give explicit estimates on the parameters $\rho_\ell$ and $r_\ell$.
\begin{Assumption} \label{ass:stronger}
Let $X$ be is an isotropic random vector in $\R^d$. Assume that there are constant $c_1$ and $C$ such that for every $t \in \R^d$ and every $1 \leq p \leq c_1\log d$,
\begin{equation} \label{eq:log-d-moments}
\|\inr{t,X}\|_{L_p} \leq C\sqrt{p}\|\inr{t,X}\|_{L_2} = C\sqrt{p}\|t\|_2~;
\end{equation}
the value of the constant $c_1$ is given in Theorem \ref{thm:expectation}.

Assume further that $\|f^*(X)-Y\|_{L_4} \leq \sigma_4$ for a known constant $\sigma_4$.
\end{Assumption}
In other words, Assumption \ref{ass:stronger} means that linear forms satisfy a sub-Gaussian moment growth, but
only up to a rather low exponent---logarithmic in the dimension of the underlying space. This moment assumption is a sufficient and almost necessary condition for the celebrated basis pursuit procedure to have
a unique minimizer (see \cite{LeMe16ab}), and as such, it is a natural assumption when studying such sparsity-driven bounds.
Note that even with Assumption \ref{ass:stronger} replacing Assumption \ref{ass:moment}, the fact that the `noise' $\xi=f^*(X)-Y$ may only be in $L_4$ means that there is not hope that
$$
\sup_{t \in T} \left|\frac{1}{\sqrt{N}} \sum_{i=1}^N \xi_i \inr{t,X_i} - \E \xi \inr{t,X} \right|~,
$$
exhibits a fast tail decay, even in the extreme case when $|T|=1$. This indicates why regularized risk minimization can only perform with a rather weak accuracy/confidence tradeoff in such situations.

On the other hand, \eqref{eq:log-d-moments} suffices to obtain bounds on the \emph{expectation} of empirical and multiplier processes, as long as the indexing set has enough symmetries, and a suitable bound on the expectation suffices for the analysis of regularized tournaments.

\begin{Definition} \label{def:K-unconditional}
Given a vector $x=(x_i)_{i=1}^n$, let $(x_i^*)_{i=1}^n$ be the non-increasing rearrangement of $(|x_i|)_{i=1}^n$.

The normed space $(\R^d,\| \ \|)$ is $K$-unconditional with respect to the basis $\{e_1,\ldots,e_d\}$ if for every $x \in \R^d$ and every permutation of $\{1,\ldots,n\}$,
$$
\left\|\sum_{i=1}^d x_i e_i\right\| \leq K \left\|\sum_{i=1}^d x_{\pi(i)}e_i\right\|,
$$
and if $y \in \R^d$ and $x_i^* \leq y_i^*$ for $1 \leq i \leq d$ then
$$
\left\|\sum_{i=1}^d x_i e_i\right\| \leq K \left\|\sum_{i=1}^d y_ie_i\right\|~.
$$
\end{Definition}

There are many natural examples of $K$-unconditional spaces, most notably, all the $\ell_p$ spaces. Moreover,  if $(v_i^*)_{i=1}^d$ denotes the nonincreasing rearrangement of $(|v_i|)_{i=1}^d$, then the norm $\|z\|=\sup_{v \in V} \sum_{i=1}^n v_i^* z_i^*$ is $1$-unconditional. In fact, if $V \subset \R^d$ is closed under coordinate permutations and reflections (sign-changes), then $\|\cdot\| = \sup_{v \in V} |\inr{\cdot,v}|$ is $1$-unconditional in the sense of Definition \ref{def:K-unconditional}.

The following fact has been established in \cite{Men16c}:

\begin{Theorem} \label{thm:expectation}
There exists an absolute constant $c_1$ and for $K \geq 1$, $L \geq 1$ and $q_0 >2$ there exists a constant $c_2$ that depends only on $K$, $L$ and $q_0$ for which the following holds. Consider
\begin{description}
\item{$\bullet$} $V \subset \R^d$ for which the norm $\|\cdot\|=\sup_{v \in V} |\inr{v,\cdot}|$ is $K$-unconditional with respect to the basis $\{e_1,\ldots,e_d\}$,
\item{$\bullet$} $\xi \in L_{q_0}$ for some $q_0>2$,
\item{$\bullet$} an isotropic random vector $X \in \R^d$ that satisfies
$$
\max_{1 \leq j \leq d} \sup_{1 \leq p \leq c_1\log d} \frac{\|\inr{X,e_j}\|_{L_p}}{\sqrt{p}} \leq L~.
$$
\end{description}
If $(X_i,\xi_i)_{i=1}^N$ are independent copies of $(X,\xi)$ then
$$
\E \sup_{v \in V} \left|\frac{1}{\sqrt{N}} \sum_{i=1}^N \left(\xi_i \inr{X_i,v} - \E \xi \inr{X,v}\right) \right| \leq c_2 \|\xi\|_{L_{q_0}} \ell_*(V)~,
$$
where $\ell_*(V)=\E \sup_{v \in V} \sum_{i=1}^d g_i v_i$ and $G=(g_i)_{i=1}^d$ is a standard Gaussian vector in $\R^d$.
\end{Theorem}

Therefore, as long as $V$ is sufficiently symmetric and linear forms
exhibit a sub-Gaussian moment growth up to $p \sim \log d$, the
expectations of empirical and multiplier processes indexed by $V$
behave as if $X$ were the standard Gaussian vector and $\xi$ were
independent of $X$. In the cases we are interested in the indexing
sets have enough symmetries, and since $\xi \in L_4$, the conditions
of Theorem \ref{thm:expectation} hold for $q_0=4$.

\subsection{The tournament {\sc lasso}}
\label{sec:lassoproof}

In this section we prove Theorem \ref{thm:LASSO-intro}, the
performance bound of the ``tournament {\sc lasso}'' procedure.


The proof of Theorem \ref{thm:LASSO-intro} follows from Theorem \ref{thm:main}, combined with explicit estimates on the parameters $\rho_\ell$, $r_{\ell,1}$ and $r_{\ell,2}$, as we now show.

Note that for any $h=\inr{t_0,\cdot}$ and every $\rho,r>0$,
$$
\F_{h,\rho,r} = \{\inr{t,\cdot} \in \R^d : t \in \rho B_1^d \cap r B_2^d\}~,
$$
where $B_p^d$ is the unit ball of the normed space $(\R^d, \| \ \|_p)$.

Hence, for a fixed radius $\rho$ the parameters $r_E$ and $\ol{r}_{\M}$ are defined using the fixed-point conditions
\begin{equation} \label{eq:LASSO-r-E}
\E \sup_{t \in \rho B_1^d \cap r B_2^d} \left| \frac{1}{\sqrt{N}} \sum_{i=1}^N \eps_i \inr{t,X_i} \right|  \leq \kappa \sqrt{N} r
\end{equation}
and
\begin{equation} \label{eq:LASSO-r-M}
\E \sup_{t \in \rho B_1^d \cap r B_2^d} \left| \frac{1}{\sqrt{N}} \sum_{i=1}^N \eps_i \xi_i \inr{t,X_i} \right|  \leq \kappa \sqrt{N} r^2
\end{equation}
respectively. The indexing set $V_{\rho,r}=\rho B_1^d \cap r B_2^d$ is invariant under coordinate permutations and sign reflections, and therefore satisfies the conditions of Theorem \ref{thm:expectation}. Hence, an upper bound on $r_E$ follows if
\begin{equation} \label{eq:gaussian-r-E}
\E \sup_{t \in \rho B_1^d \cap r B_2^d } \sum_{i=1}^N g_i t_i  \leq \kappa \sqrt{N} r~,
\end{equation}
while for an upper estimate on $\ol{r}_{\M}$ it suffices to ensure that
\begin{equation} \label{eq:gaussian-r-M}
\|\xi\|_{L_4} \E \sup_{t \in \rho B_1^d \cap r B_2^d } \sum_{i=1}^N g_i t_i  \leq \kappa \sqrt{N} r^2~.
\end{equation}
Equations \eqref{eq:gaussian-r-E} and \eqref{eq:gaussian-r-M} cannot be improved; they are tight bounds on \eqref{eq:LASSO-r-E} and \eqref{eq:LASSO-r-M} when, for example, $X=(g_1,\ldots,g_d)$ and $\xi$ is a Gaussian variable that is independent of $X$.

The added value in \eqref{eq:gaussian-r-E} and \eqref{eq:gaussian-r-M} is that if $r$ satisfies these inequalities then, necessarily, $\max\{\lambda_{\Q},\lambda_{\M}\} \leq r$ (for a well chosen constant $\kappa$). This is an immediate consequence of Sudakov's inequality, which implies that for some absolute constant $c>0$, for any $T \subset \R^d$ and any
$\eps>0$,
$$
\eps \sqrt{\log {\cal M}(T,\eps B_2^d)} \leq c \E \sup_{t \in T} \sum_{i=1}^d g_i t_i \equiv \ell_*(T)~.
$$
Thus, when applied to the definition on $\lambda_{\Q}$ one obtains
$$
\log {\cal M}(\rho B_1^d \cap r B_2^d, \eta r B_2^d) \leq \frac{\ell_*^2(\rho B_1^d \cap r B_2^d)}{(\eta r)^2} \leq \kappa^2 N~,
$$
that is, it suffices that
$$
\ell_*(\rho B_1^d \cap r B_2^d) \leq \eta \kappa \sqrt{N} r~,
$$
which is precisely the type of condition in \eqref{eq:gaussian-r-E}.

With those estimates in hand, we are now able to explain the choice of parameters $\rho_\ell$, $r_{\ell,1}$, $r_{\ell,2}$ and $\lambda_\ell$ for the tournament {\sc lasso}. This requires several observations that have been established in \cite{LeMe16a}.

First, the requirement that $\Delta(\rho,r) \geq 4\rho/5$ forces some constraint on the choice of $\rho$ and $r$. To simplify things, assume that $t_0={\rm argmin}_{t \in \R^d} \E(\inr{X,t}-Y)^2$ is supported on $I \subset \{1,\ldots,d\}$ and that $|I| \leq s$. Recall that by the definition of $\Delta(\rho,r)$, the fact that $\Psi(t)=\|t\|_1$ and since $X$ is isotropic, it suffices to consider vectors $t \in \R^d$ for which $\|t-t_0\|_1=\rho$ and $\|t-t_0\|_2 \leq r$.
For such $t$,
$$
\|t\|_1 - \|t_0\|_1 = \sum_{i \in I^c} |t_i| + \sum_{i \in I} \left(|t_i| - |(t_0)_i|\right) \geq \sum_{i \in I^c} |t_i| - \sum_{i \in I} |t_i-(t_0)_i|~,
$$
and since $|I| \leq s$,
$$
\sum_{i \in I} |t_i-(t_0)_i| \leq \sqrt{|I|}\|t-t_0\|_2 \leq \sqrt{s}r~.
$$
Therefore,
$$
\sum_{i \in I^c} |t_i| = \sum_{i \in I^c} |t_i-(t_0)_i| = \sum_{i=1}^n |t_i-(t_0)_i| - \sum_{i \in I} |t_i-(t_0)_i| \geq \rho - \sqrt{s}r~.
$$
On the other hand, there is a functional $z$ that is norming for both $t_0$ and $P_{I^c}t =\sum_{i \in I^c} t_i e_i$; hence,
\begin{eqnarray*}
z(t-t_0) & \geq & z(P_{I^c}(t-t_0)) - \sum_{i \in I} |t_i-(t_0)_i| = \sum_{i \in I^c}|t_i-(t_0)_i| - \sum_{i \in I} |t_i-(t_0)_i|
\\
& \geq & \rho - 2\sqrt{s}r \geq \frac{4\rho}{5}
\end{eqnarray*}
as long as $s \lesssim (\rho/r)^2$.

This shows that, as long as the ratio $\rho/r$ is larger than the square-root of the degree of sparsity of vectors we are interested in, $\Delta(\rho,r) \geq (4/5)\rho$ as our procedure requires. A similar observation is true if $t_0$ is not sparse, but rather well approximated by an $s$-sparse vector (see \cite{LeMe16a} for a detailed argument).

Set $k=(\rho/r)^2$ and assume without loss of generality that $k$ is an integer. We also restrict ourselves to values $1 \leq k \leq d$, intuitively because the above implies that $(\rho/r)^2$ should capture the degree of sparsity. Recall that
$$
\ell_*(\rho B_1^d \cap r B_2^d) = r \ell_*(\sqrt{k}B_1^d \cap B_2^d) \leq Cr \sqrt{k \log (ed/k)} = C\rho \sqrt{\log(edr^2/\rho^2)}
$$
(see, e.g. \cite{LeMe16a} for the standard proof). Hence,  \eqref{eq:gaussian-r-E} becomes
\begin{equation} \label{eq:gaussian-r-E-1}
C\rho \sqrt{\log(edr^2/\rho^2)}  \leq \kappa \sqrt{N} r~,
\end{equation}
while \eqref{eq:gaussian-r-M} implies
\begin{equation} \label{eq:gaussian-r-M-1}
\|\xi\|_{L_4} \cdot C\rho \sqrt{\log(edr^2/\rho^2)}  \leq \kappa \sqrt{N} r^2.
\end{equation}
We consider only the case $N \leq C d$, which is the more interesting range in sparse recovery---when the number of given linear measurements is significantly smaller than the dimension of the underlying space. An argument following the same path may be used when $N \geq C d$ and we omit it.

It follows from a rather tedious computation that \eqref{eq:gaussian-r-E-1} holds provided that
\begin{equation} \label{eq:gaussian-r-E-2}
r \geq c\frac{\rho}{\kappa \sqrt{N}} \sqrt{\log\left(\frac{cd}{\kappa N}\right)}~,
\end{equation}
and it follows from \eqref{eq:gaussian-r-M-1} that
\begin{equation} \label{eq:gaussian-r-M-2}
r^2 \geq c\rho \frac{\|\xi\|_{L_4}}{\sqrt{N}} \sqrt{\log\left(c\frac{\|\xi\|_{L_4} d}{\sqrt{N} \rho}\right)}
\end{equation}
as long as $\|\xi\|_{L_4} d /\sqrt{N} \rho \geq c^\prime$.

Using the constraint that $\rho/r \geq c\sqrt{s}$, it is evident
from \eqref{eq:gaussian-r-E-2} that

$s \leq c(L) N/\log\left({ed}/{N}\right)$, and that
$$
\frac{1}{s} \geq c_1(L) \frac{r^2}{\rho^2} \gtrsim \frac{1}{\rho} \cdot \frac{\|\xi\|_{L_4}}{\sqrt{N}} \sqrt{\log\left(e\frac{\|\xi\|_{L_4} d}{\sqrt{N} \rho}\right)}~.
$$
Therefore, to have a `legal' choice of $\rho$ and $r$, we must have
$$
N \geq c_2(L) s \log \left(\frac{ed}{s}\right)~,
$$
and
$$
\rho \geq c_3(L) \frac{s}{\sqrt{N}} \|\xi\|_{L_4} \cdot \sqrt{\log\left(\frac{ed}{s}\right)}~.
$$
This naturally leads to the choices made in Section \ref{sec-lasso-intro}: set
$$
\F_\ell = \{ t : \exists v, \ |{\rm supp}(v)| \leq d/2^{\ell-1}, \ \|t-v\|_1 \leq \rho_\ell\}
$$
to be the set of vectors that are `well-approximated' by $d/2^{\ell-1}$ sparse vectors. For every $\ell$ let $s = d/2^{\ell-1}$; if $s \geq c(L) N/\log\left(\frac{ed}{N}\right)$, set $\rho_\ell=r_\ell=\infty$. If the reverse inequality holds, set
$$
\rho_\ell = c(L) \frac{d}{2^\ell \sqrt{N}} \|\xi\|_{L_4}\sqrt{\log\left(e2^\ell \right)} \sim_L \frac{d\sqrt{\ell}}{2^\ell \sqrt{N}} \|\xi\|_{L_4}~,
$$
and for that choice $\rho_\ell$, the required value of $r_\ell$ is
$$
r_\ell \geq c(L)\|\xi\|_{L_4} \sqrt{\frac{s}{N} \log\left(\frac{ed}{s}\right)} \sim_L \|\xi\|_{L_4} \sqrt{\frac{d \ell }{2^\ell N}}~.
$$
Finally, let $\wh{r}_\ell \geq r_\ell$ and recall that $\lambda_\ell \sim_L \wh{r}_\ell^2/\rho_\ell$.  Applying Theorem \ref{thm:main}, these choices complete the proof of Theorem \ref{thm:LASSO-intro}.
\endproof

\subsection{The tournament {\sc slope}}

As a second example, we present and analyze a ``tournament'' version of
the regularized risk minimization procedure {\sc slope}.
{\sc slope} is defined using a set of non-increasing weights $(\beta_i)_{i=1}^d$. The corresponding norm is
$$
\Psi(z)=\sum_{i=1}^d \beta_i z_i^*~,
$$
where as always, $(z_i^*)_{i=1}^d$ denotes the non-increasing rearrangement of $(|z_i|)_{i=1}^d$. Clearly, {\sc slope} is a generalized version of {\sc lasso}, as the latter is given by the choice $\beta_i=1$ for $1
\leq i \leq d$.

Just like the {\sc lasso}, most of the known results on the performance of {\sc slope} hold only when both the random vector $X$ and the target $Y$ have well behaved tails.

The tournament {\sc slope} we present below is defined for the penalty $\Psi(z)=\sum_{i=1}^d \beta_i z_i^*$, where $\beta_i \leq C\sqrt{\log(ed/i)}$. We obtain the following
performance bound:

\begin{Theorem} \label{thm:tournament-slope}
For $L \geq 1$ there exist constants $c_0,\ldots,c_5$ that depend only on $L$ and for which the following holds.
Let $X$ satisfy Assumption \ref{ass:stronger}, set $t_0={\rm argmin}_{t \in \R^d} \E(Y-\inr{t,X})^2$ and assume that $\|Y-\inr{t_0,X}\|_{L_4} \leq \sigma$. Assume further that there is $v$ that is $s$-sparse such that
$$
\|t_0-v\|_1 \leq c_1(L) \sigma \cdot \frac{s\log(ed/s)}{\sqrt{N}}~.
$$
If  $N \geq c_2(L)s \log(ed/s)$ and
$$
\wh{r} \geq c_3(L) \sigma \sqrt{\frac{s}{N}\log\left(\frac{ed}{s}\right)}~,
$$
then with probability at least
$$
1-2\exp\left(-c_4(L)N \min\left\{1,\left(\frac{\wh{r}}{\sigma}\right)^2\right\}\right)~,
$$
the tournament {\sc slope} produces $\wh{t}$ that satisfies
$$
\|\wh{t}-t_0\|_2 \leq c_5(L) \sigma \sqrt{\frac{s}{N}\log\left(\frac{ed}{s}\right)},
$$
$$
\Psi(\wh{t}-t_0)=\sum_{i=1}^d (\wh{t}-t_0)_i^* \sqrt{\log(ed/i)} \leq c_5(L) \sigma \frac{s}{\sqrt{N}} \log \left(\frac{ed}{s}\right)~,
$$
and
$$
\E \left( \left(\inr{\wh{t},X}-Y\right)^2 | (X_i,Y_i)_{i=1}^N \right) \leq \E \left(\inr{t_0,X}-Y\right)^2 + c_5(L)\wh{r}^2~.
$$
\end{Theorem}
The estimate corresponds to the optimal accuracy/confidence tradeoff any procedure can attain even if the learner knows that $t_0$ is $s$-sparse. Moreover, in the heavy-tailed situations we study here, the performance of {\sc slope} is significantly weaker than in Theorem \ref{thm:tournament-slope}.

The argument we use here is similar to the one used for the tournament {\sc lasso}, and so we skip most of the details.

In tournament {\sc slope} one selects $\beta_i \leq
c_0\sqrt{\log(ed/i)}$ and therefore the corresponding indexing set is contained in $$
V_{\rho,r}=\rho {\cal B} \cap r B_2^d = \left\{ v \in \R^d : \|v\|_2 \leq r \ {\rm and} \ \sum_{i=1}^d v_i^* \sqrt{\log(ed/i)} \leq \rho/c_0 \right\}~.
$$

Because $V_{\rho,r}$ has enough symmetries, one may apply Theorem \ref{thm:expectation}, leading to an upper bound on $r_E$ when
\begin{equation} \label{eq:gaussian-slope-E-1}
\E \sup_{v \in V_{\rho,r}} \sum_{i=1}^d g_i v_i \leq \kappa \sqrt{N}r~.
\end{equation}
Also, to estimate $r_{\M}$ it suffices to ensure that
\begin{equation} \label{eq:gaussian-slope-M-1}
\|\xi\|_{L_4} \E \sup_{v \in V_{\rho,r}} \sum_{i=1}^d g_i v_i \leq \kappa \sqrt{N}r^2~.
\end{equation}
Next, one may verify (see Lemma 4.3 in \cite{LeMe16a}) that if we set $B_s = \sum_{i \leq s} \beta_i/\sqrt{i}$ and if $B_s \lesssim r/\rho$, then $\Delta(\rho,r) \geq (4/5)\rho$ for centres that are `well approximated' by $s$-sparse vectors. Also, for our choice of $\beta_i$, $B_s \lesssim C \sqrt{s \log(ed/s)}$. Hence, for a fixed degree of sparsity $1 \leq s \leq d$, one has the constraint that
\begin{equation} \label{eq:sparse-cond-SLOPE}
\frac{r}{\rho} \geq C_1 \sqrt{s \log(ed/s)}
\end{equation}
for a constant $C_1$ that depends only on $c_0$.

Following the same path used for the tournament {\sc lasso} let
$$
\F_\ell = \{ t : \exists v, \ |{\rm supp}(v)| \leq d/2^{\ell-1}, \ \Psi(t-v) \leq \rho_\ell\}.
$$
For every $1 \leq \ell \leq \log_2 K$ let $s=d/2^{\ell-1}$. There is a nontrivial choice of $\rho$ and $r$ only when $s \lesssim_L N/\log(ed/N)$; otherwise, $\rho=r=\infty$ as one would expect. When $s \lesssim_L N/\log(ed/N)$, we follow the computation in \cite{LeMe16a} and set
$$
\rho_\ell \sim_L \|\xi\|_{L_4} \frac{s}{\sqrt{N}} \log \left(\frac{ed}{s}\right) \sim_L \|\xi\|_{L_4} \frac{d \ell}{2^\ell\sqrt{N}}~,
$$
and
$$
r_{\ell} \sim_L \|\xi\|_{L_4} \sqrt{\frac{s}{N}\log\left(\frac{ed}{s}\right)}~.
$$
Finally, fix $\wh{r}_\ell \geq r_\ell$ and set $r_{\ell,1},r_{\ell,2} \sim_L \wh{r}_\ell$, $\lambda_\ell \sim_L \wh{r}_\ell^2/\rho$. Applying Theorem \ref{thm:main} for these choices completes the proof of Theorem \ref{thm:tournament-slope}.
\endproof

\appendix
\section{Additional proofs}
%
%

The proofs of Lemma \ref{lemma:basic-combining-loss-and-reg} and Lemma \ref{lemma:components-of-main} are, in fact, the same as in \cite{LeMe16a} and \cite{LuMe16}, respectively. The minor modifications to the original proofs are presented in this appendix solely for the sake of completeness and not in full detail.
\subsection*{Proof of Lemma \ref{lemma:basic-combining-loss-and-reg}}
The proof of Lemma \ref{lemma:basic-combining-loss-and-reg} follows the same path as that of Theorem 3.2 in \cite{LeMe16a}. Let us begin by examining
$$
(*) = \Q_{f,f^*}(j)+\M_{f,f^*}(j)+\lambda(\Psi(f)-\Psi(f^*))
$$
in the set $\{f \in \F : \Psi(f-f^*)=\rho\}$. If $\Psi(f-f^*)=\rho$  one should consider two cases. First, if $\|f-f^*\|_{L_2} \geq r$ then by the triangle inequality for $\Psi$, and since $\Q_{f,f^*}(j) \geq C\|f-f^*\|_{L_2}^2$ and $\M_{f,f^*}(j) \geq -(C/4)\|f-f^*\|_{L_2}^2$, we have
\begin{eqnarray} \label{eq:large-L-2}
(*)
& \geq & C\|f-f^*\|_{L_2}^2 - \frac{C}{4}\|f-f^*\|_{L_2}^2 - \lambda \Psi(f-f^*)
\\
& \geq & \frac{3C}{4} \|f-f^*\|_{L_2}^2 -\lambda \rho \geq \frac{3C}{4}r^2 - \lambda \rho > 0~, \nonumber
\end{eqnarray}
provided that
\begin{equation} \label{eq:lambda-upper}
\lambda \leq \frac{3C}{4} \cdot \frac{r^2}{\rho}~.
\end{equation}
If, on the other hand, $\|f-f^*\|_{L_2} \leq r$, then $\Q_{f,f^*}(j) \geq 0$ and $\M_{f,f^*}(j) \geq -(C/4)r^2$. Therefore,
$$
(*)
\geq -\frac{C}{4}r^2 + \lambda(\Psi(f)-\Psi(f^*))~.
$$
Fix $v \in {\cal B}_{f^*}(\rho/20)$ and write $f^*=u+v$; thus $\Psi(u) \leq \rho/20$.  Set $z$ to be a linear functional that is norming for $v$ and observe that for any $f \in E$,
\begin{eqnarray} \label{eq:small-L-2a}
\Psi(f)-\Psi(f^*) & \geq & \Psi(f) - \Psi(v) - \Psi(u) \geq z(f-v)-\Psi(u) \geq z(f-f^*)-2\Psi(u) \nonumber
\\
& \geq & z(f-f^*) - \frac{\rho}{10}~.
\end{eqnarray}
Hence, if $f^* \in \F_\ell$ and $f \in \F \cap {\cal B}_{f^*}(\rho) \cap D_{f^*}(r)$ then optimizing the choices of $v$ and of $z$, $z(f-f^*) \geq \Delta_\ell(\rho,r)$; thus
\begin{equation} \label{eq:small-L-2b}
\Psi(f)-\Psi(f^*) \geq \Delta_\ell(\rho,r) -\frac{\rho}{10} \geq \frac{7}{10} \rho~.
\end{equation}
And, if
\begin{equation} \label{eq:lambda-lower}
\lambda \geq \frac{C}{2} \cdot \frac{r^2}{\rho}~,
\end{equation}
we have that
$$
(*) \geq -\frac{C}{4}r^2 + \lambda \cdot \frac{7}{10} \rho > 0~.
$$
In other words, if $\lambda$ is chosen to satisfy both \eqref{eq:lambda-upper} and \eqref{eq:lambda-lower}, $f \in \F$ and $\Psi(f-f^*)=\rho$, it follows that
$$
\Q_{f,f^*}(j) + \M_{f,f^*}(j) + \lambda (\Psi(f)-\Psi(f^*))>0~.
$$
Next, if $\Psi(f-f^*) > \rho$, there are $\theta \in (0,1)$ and $h \in \F$ that satisfy
$$
\Psi(h-f^*) = \rho \ \ {\rm and} \ \ \theta (f-f^*) = h-f^*.
$$
If $\|h-f^*\|_{L_2} \geq r$, then by the triangle inequality for $\Psi$ followed by \eqref{eq:large-L-2},
\begin{eqnarray*}
(*) & \geq & \frac{1}{\theta^2} \Q_{h,f^*}(j) + \frac{1}{\theta}\left(\M_{h,f^*}(j) - \lambda \Psi(h-f^*) \right)
\\
& \geq & \frac{1}{\theta} \left(\Q_{h,f^*}(j)+\M_{h,f^*}(j) - \lambda \Psi(h-f^*) \right) > 0~.
\end{eqnarray*}
If, on the other hand, $\|h-f^*\|_{L_2} \leq r$, then
\begin{eqnarray*}
(*)
& \geq & \frac{1}{\theta} \M_{h,f^*}(j) + \lambda(z(f-f^*)-2\Psi(u))
\\
& \geq & \frac{1}{\theta} \left(\M_{h,f^*}(j) +
  \lambda\left(z(h-f^*)-2\theta \Psi(u)\right)\right)
\\
& \geq & \frac{1}{\theta} \left(\M_{h,f^*}(j) + \lambda\left(z(h-f^*)-2\Psi(u)\right)\right)>0~,
\end{eqnarray*}
because $0 \leq \theta <1$ and using \eqref{eq:small-L-2a}.

Now, all that remains is to control $f \in \F \cap {\cal B}_{f^*}(\rho)$ and show that if $\|f-f^*\|_{L_2} \geq r$, then
$$
\Q_{f,f^*}(j)+\M_{f,f^*}(j)+\lambda(\Psi(f)-\Psi(f^*))>0~.
$$
This follows from \eqref{eq:large-L-2}.
\endproof

\subsection*{Proof of Lemma \ref{lemma:components-of-main}}
The first part of Lemma \ref{lemma:components-of-main} is identical to Lemma 5.1 from \cite{LuMe16}, with the trivial modification that the constant $-C/4$ replaces $-3C/4$ used in \cite{LuMe16}. The second part of Lemma \ref{lemma:components-of-main} was not needed in \cite{LuMe16}, but its proof follows the same path as Lemma 5.1 from \cite{LuMe16}.

Set $r=\wh{r}_\ell$ and fix $f \in \F$ that satisfies $\|f-f^*\|_{L_2} \leq r$. Recall that $m=N/n$ and that $\sqrt{n/N} \leq \sqrt{\theta}r/\sigma_4$ for a well-chosen constant $\theta$ that depends only on $L$ and $\tau$. Set $U=(f-f^*)(X) \cdot (f^*(X)-Y)$ and observe that
$$
\M_{f,f^*}=\frac{1}{m}\sum_{i=1}^m U_i.
$$
It follows from the convexity of $\F$ that $\E U = \E (f-f^*)(X) \cdot (f^*(X)-Y) \geq 0$; therefore, $\M_{f,f^*} \geq \M_{f,f^*}-\E \M_{f,f^*}$. Also,
$$
Pr \left( \left| \M_{f,f^*} - \E \M_{f,f^*} \right| > t\right) \leq t^{-1}\E | \M_{f,f^*} - \E \M_{f,f^*}|,
$$
and by a straightforward symmetrization argument,
$$
\E | \M_{f,f^*} - \E \M_{f,f^*}| \leq 2\E \left|\frac{1}{m}\sum_{i=1}^m \eps_i U_i \right| \leq \frac{2}{\sqrt{m}} (\E |U|^2)^{1/2}.
$$
Applying Assumption \ref{ass:moment}, it is evident that
$$
(\E|U|^2)^{1/2} \leq \|f^*(X)-Y\|_{L_4} \cdot \|f-f^*\|_{L_4} \leq L r \sigma_4,
$$
and thus
$$
Pr \left( \left| \M_{f,f^*} - \E \M_{f,f^*} \right| > t\right) \leq \frac{2Lr \sigma_4}{t \sqrt{m}} =  \frac{2L \sigma_4 r \sqrt{n}}{t \sqrt{N}} \leq \frac{\tau}{3},
$$
where we use the fact that $\sqrt{n/N} \leq \sqrt{\theta} r/\sigma_4$ and select $t = Cr^2/8$ and $\theta=\theta(\tau,L)$. Therefore,
$$
Pr\left(\M_{f,f^*} \leq - (C/8)r^2\right) \leq \frac{\tau}{3},
$$
and with probability at least $1-2\exp(-c\tau^2n)$,
\begin{equation} \label{eq:M-single}
\left|\left\{j : \M_{f,f^*}(j) \geq -(C/8)r^2 \right\}\right| \geq (1-\tau/2)n.
\end{equation}
The rest of the argument is identical to the proof of Lemma 5.1 from \cite{LuMe16}: let $\cH$ be a maximal separated subset of $\F \cap {\cal B}_{f^*}(\rho) \cap D_{f^*}(r)$ with respect to the $L_2$ norm, of cardinality $\exp(c\tau^2n/2)$, and with the following property: for any $f \in \F \cap {\cal B}_{f^*}(\rho) \cap D_{f^*}(r)$ there is $h \in \cH$ for which
\begin{equation} \label{eq:M-approx}
\|f-h\|_{L_2} \leq \eps \ \ {\rm and} \ \ \E(f^*(X)-Y)(f(X)-h(X)) \geq 0;
\end{equation}
here $\eps$ denotes the mesh of the net. The existence of such a separated set is established in \cite{LuMe16} (see Lemma 5.3), and one may show that the mesh $\eps$ is a small proportion of $r$.

By \eqref{eq:M-single}, we have that with probability at least $1-2\exp(-c\tau^2n/2)$, for every $h \in \cH$
\begin{equation} \label{eq:M-single-net}
\left|\left\{j : \M_{h,f^*}(j) \geq -(C/8)r^2 \right\}\right| \geq (1-\tau/2)n.
\end{equation}

For every $f \in \F \cap {\cal B}_{f^*}(\rho) \cap D_{f^*}(r)$ let $\pi f \in \cH$ be as in \eqref{eq:M-approx}, and at the heart of  the proof of Lemma 5.4 in \cite{LuMe16} is that with probability at least $1-2\exp(-c_1\tau^2 n)$,
\begin{equation} \label{eq:M-uniform}
\sup_{f \in \F \cap {\cal B}_{f^*}(\rho) \cap D_{f^*}(r)} \left|\left\{j : \M_{f,f^*}(j) - \M_{\pi f,f^*}(j) \leq -(C/8)r^2 \right\}\right| \leq \frac{\tau n}{2}.
\end{equation}
Combining \eqref{eq:M-single-net} and \eqref{eq:M-uniform}, there is an event of probability at least $1-2\exp(-c_2\tau^2n)$ on which for any $f \in \F \cap {\cal B}_{f^*}(\rho) \cap D_{f^*}(r)$ there is a set of coordinate blocks $(I_j)_{j \in J}$, of cardinality $|J| \geq (1-\tau)n$ and for $j \in J$,
$$
\M_{f,f^*}(j) \geq \M_{\pi f , f^*}(j) + \left(\M_{f,f^*}(j) - \M_{\pi f ,f ^*}(j)\right) \geq -\frac{C}{4}r^2.
$$
\endproof

\subsection*{Acknowledgements}
We thank the referees for valuable suggestions that significantly helped
us improve the presentation.

\bibliographystyle{plain}

\end{document}